\documentclass[12pt]{article}                
\usepackage{graphicx}
\usepackage{psfrag}
\usepackage{amsmath, amssymb}
\usepackage{mathptmx}
\usepackage{ntheorem}

\DeclareMathOperator{\tr}{tr}
\DeclareMathOperator{\I}{I}
\DeclareMathOperator{\diag}{diag}
\DeclareMathOperator{\Cov}{Cov}

\newcommand{\A}{{\cal A}}

\newcommand{\N}{{\cal N}}

\newcommand{\DD}{{\cal D}}
\newcommand{\PP}{{\cal P}}

\newcommand{\eps}{{\varepsilon}}

\newcommand{\smsp}{\hspace{0.3mm}}
\newcommand{\e}{\mathbb{E}}

\newcommand{\p}{\mathbb{P}}

\newcommand{\Reals}{\mathbb{R}}
\newcommand{\Natural}{\mathbb{N}}
\newcommand{\la}{\langle}
\newcommand{\ra}{\rangle}
\newcommand\qed{\hfill\hbox{\rlap{$\sqcap$}$\sqcup$}}

\newcommand{\bigO}[1]{\ensuremath{\mathop{}\mathopen{}\mathcal{O}\mathopen{}\left(#1\right)}}

\newtheorem{lemma}{Lemma}
\newtheorem{theorem}{Theorem}

\theoremstyle{nonumberplain}
\newcommand\specialref{}

\begin{document}

\title{Free energy in the Potts spin glass}
\author{Dmitry Panchenko\thanks{\textsc{\tiny Department of Mathematics, University of Toronto, panchenk@math.toronto.edu. Partially supported by NSERC.}}\\
}
\date{}
\maketitle
\begin{abstract}
We study the Potts spin glass model, which generalizes the Sherrington-Kirkpatrick model to the case when spins take more than two values but their interactions are counted only if the spins are equal. We obtain the analogue of the Parisi variational formula for the free energy, with the order parameter now given by a monotone path in the set of positive-semidefinite matrices. The main idea of the paper is a novel synchronization mechanism for blocks of overlaps. This mechanism can be used to solve a more general version of the Sherrington-Kirkpatrick model with vector spins interacting through their scalar product, which includes the Potts spin glass as a special case. As another example of application, one can show that Talagrand's bound for multiple copies of the mixed $p$-spin model with constrained overlaps is asymptotically sharp. We will consider these problems in the subsequent paper and illustrate the main new idea on the technically more transparent case of the Potts spin glass.
\end{abstract} 
\vspace{0.5cm}

\section{Introduction and main results}

The Hamiltonian of the classical Sherrington-Kirkpatrick model introduced in \cite{SK} is a random function of the $N\geq 1$ spins taking values $\pm 1,$
\begin{equation}
\sigma = (\sigma_1,\ldots,\sigma_N) \in \bigl\{-1,+1 \bigr\}^N,
\end{equation}
given by the quadratic form
\begin{equation}
H_N(\sigma) = \frac{1}{\sqrt{N}} \sum_{1\leq i<j\leq N} g_{ij}\sigma_i \sigma_j,
\label{SKH}
\end{equation}
where the interaction parameters $(g_{ij})$ are independent standard Gaussian random variables. One common interpretation of this Hamiltonian (see e.g. \cite{SDean}) is related to the following so-called Dean's problem (another variant was named Shakespeare's problem in the classic book on spin glasses \cite{MPV}). Given a group of $N$ students, the parameter $g_{ij}$ represents how much the students $i$ and $j$ like or dislike each other and the configuration $\sigma$ represents a possible assignment of the students to the two dormitories, labeled $\pm1$, by their dean. If a pair of students $(i,j)$ is assigned to the same dormitory, their interaction $g_{ij}$ is counted in (\ref{SKH}) with the plus sign, $\sigma_i \sigma_j=+1$, otherwise, it is counted with the minus sign, $\sigma_i \sigma_j=-1$. The Hamiltonian $H_N(\sigma)$ is viewed as the global comfort function and one is then interested in understanding the behaviour of its maximum, which can be related to the problem of computing the free energy. 

There is a natural generalization of the Dean's problem to the case of $\kappa$ dormitories for $\kappa \geq 2$, called the Potts spin glass, which has been studied extensively in the physics literature (see e.g. \cite{ES, ES2, Nishimori, Gross, Santis, Calta, Marinari}), although the formula for the free energy that we prove in this paper never appeared in full generality. The spin configurations in this model are given by 
\begin{equation}
\sigma = (\sigma_1,\ldots,\sigma_N) \in  \bigl\{1,2, \ldots, \kappa\bigr\}^N,
\end{equation}
and the Hamiltonian is defined by
\begin{equation}
H_N(\sigma) = \frac{1}{\sqrt{N}} \sum_{1\leq i,j\leq N} g_{ij}\I(\sigma_i = \sigma_j).
\label{DeanH}
\end{equation}
In physics, $\kappa$ values of spins are called `orientations' or `states'. Compared to (\ref{SKH}), the interaction term $g_{ij}$ in (\ref{DeanH}) is counted with the factor $\I(\sigma_i = \sigma_j)\in\{0,1\}$ instead of $2\I(\sigma_i = \sigma_j)-1\in\{-1,+1\}$. This transformation only rescales the Hamiltonian and shifts it by a constant (random) factor, so it is irrelevant to the computation of the free energy. Also, for convenience, we sum over all pairs of indices $(i,j)$. Except for these minor differences, the Hamiltonian (\ref{DeanH}) represents the comfort function of the Dean's problem with $\kappa$ dormitories with the students still counted as `friends' or `enemies' depending on whether they are assigned to the same dormitory or not. In the standard version of the model one also consider the case where the Gaussian variables $g_{ij}$ have non-zero mean of the order $1/\sqrt{N}$ but, for simplicity of notation, we will focus only on the conceptual difficulties related to the purely random part of the disorder. Our main goal will be to find the general formula for the limit of the free energy
\begin{equation}
F_N = \frac{1}{N}\smsp \e \log \sum_{\sigma} \exp \beta H_N(\sigma)
\label{FE}
\end{equation}
for any inverse temperature parameter $\beta>0$. 

The Potts spin glass is a special case of the following version of the Sherrington-Kirkpatrick model with vector-valued spins.  The spin configurations in this model are given by 
\begin{equation}
\sigma = (\sigma_1,\ldots,\sigma_N) \in  (\Reals^{\kappa})^N,
\end{equation}
and the Hamiltonian is defined by
\begin{equation}
H_N(\sigma) = \frac{1}{\sqrt{N}} \sum_{1\leq i,j\leq N} g_{ij}\, (\sigma_i,\sigma_j),
\label{VectorSK}
\end{equation}
where $(\sigma_i, \sigma_j)$ is the scalar product of $\sigma_i$ and $\sigma_j$. If we consider a probability measure $\mu$ on a bounded subset $\Omega\subseteq \Reals^{\kappa}$ then the free energy is defined by
\begin{equation}
F_N = \frac{1}{N}\smsp \e \log \int_{\Omega^N} \exp \beta H_N(\sigma)\, d\mu^{\otimes N}.
\end{equation}
The model (\ref{DeanH}) corresponds to the case when $\mu$ is the uniform distribution on the standard basis of $\Reals^{\kappa}$. Not to hide the main idea in the technical details, we first present the case of the  Potts spin glass, and in the subsequent paper \cite{PVS} we consider the general mixed $p$-spin version of the Sherrington-Kirkpatrick model with vector spins.

Going back to (\ref{DeanH}), it is enough to compute the limit of the free energy with fixed proportions of spins in different states. Consider the set
\begin{equation}
\DD = \Bigl\{ (d_1,\ldots, d_\kappa) \mid d_1,\ldots, d_\kappa \geq 0, \sum_{k\leq \kappa} d_k = 1 \Bigr\}
\end{equation}
of possible proportions of the states. For $d\in \DD$, we consider the set of configurations 
\begin{equation}
\Sigma(d) =  \Bigl\{ \sigma \mid \sum_{i\leq N}\I(\sigma_i=k)= N d_k \,\mbox{ for all }\, k\leq \kappa \Bigr\}
\label{SigmaD}
\end{equation}
constrained by the state sizes and define the constrained free energy by
\begin{equation}
F_N(d) = \frac{1}{N}\smsp \e \log \sum_{\sigma\in \Sigma(d)} \exp \beta H_N(\sigma).
\label{FED}
\end{equation}
By the classical Gaussian concentration inequalities, $F_N$ is approximated by $\max_{d\in \DD} F_N(d)$ and most of the work will go into the computation of this constrained free energy.

Before we write down our main result, let us describe a new phenomenon that will appear in this model that will allow us to overcome the main difficulty in the computation of the free energy (for precise formulation, see Theorem \ref{Th2} below).  As in the Sherrington-Kirkpatrick model, a crucial role will be played by the distribution of the overlaps between i.i.d. replicas $(\sigma^\ell)_{\ell\geq 1}$ from the Gibbs measure (with configurations restricted by the state sizes), only now, for a pair of replicas $\sigma^\ell$ and $\sigma^{\ell'}$, we will need to consider a $\kappa\times\kappa$ matrix of different types of overlaps
\begin{equation}
R_{\ell,\ell'}^{k,k'} = \frac{1}{N}\sum_{i\leq N}\I(\sigma_i^\ell = k) \I(\sigma_i^{\ell'}=k')
\label{Rllkk}
\end{equation} 
indexed by $k,k'\leq \kappa.$ The main novelty of the paper will be a mechanism that, by way of a small perturbation of the Hamiltonian, will force the overlap matrix 
\begin{equation}
R_{\ell,\ell'}=(R_{\ell,\ell'}^{k,k'})_{k,k'\leq \kappa}
\label{Rll}
\end{equation}
to concentrate in the infinite-volume limit $N\to\infty$ on the set of Gram matrices
\begin{equation}
\Gamma_\kappa = \Bigl\{\gamma \mid \gamma \mbox{ is a $\kappa\times\kappa$ symmetric positive-semidefinite matrix}\Bigr\}. 
\label{GammaK}
\end{equation}
Moreover, asymptotically, the entire (random) matrix $R_{\ell,\ell'}$ will become a deterministic function of its trace $\tr(R_{\ell,\ell'}),$ and this function will be non-decreasing in the sense of matrix comparison. As a result, for a model constrained to the set (\ref{SigmaD}), we will be able to describe the distribution of $R_{\ell,\ell'}$ by a sort of `quantile transformation' belonging to the set
\begin{equation}
\Pi_d = \Bigl\{ \pi\in \Pi \mid \pi(0)=0 \mbox{ and }\, \pi(1) = \diag(d_1,\ldots,d_\kappa)\Bigr\}
\label{Pid}
\end{equation}
where $\Pi$ is the space of left-continuous monotone functions (paths) in $\Gamma_\kappa$,
\begin{equation}
\Pi = \Bigl\{ \pi \colon [0,1]\to \Gamma_\kappa \mid  \pi \mbox{ is left-continuous, } \pi(x)\leq \pi(x') \mbox{ for } x\leq x' \Bigr\}.
\end{equation}
Of course, $\pi(x)\leq \pi(x')$ means that $\pi(x')-\pi(x)\in \Gamma_\kappa.$ Combined with the fact that the array $(\tr(R_{\ell,\ell'}))_{\ell,\ell'\geq 1}$ will be ultrametric by the main result in \cite{PUltra} and generated by the Ruelle probability cascades \cite{Ruelle}, this will allow us to encode the distribution of the entire array 
\begin{equation}
R= (R_{\ell,\ell'}^{k,k'})_{\ell,\ell'\geq 1, k,k'\leq \kappa}
\end{equation}
in terms of one element $\pi\in \Pi_d$ which plays the role analogous to the Parisi functional order parameter in the Sherrington-Kirkpatrick model. Notice that a priori the matrix $R_{\ell,\ell'}$ in (\ref{Rll}) is not even symmetric and, at first look, other properties mentioned above seem even less plausible. Nevertheless, a novel matrix version of the synchronization mechanism discovered in \cite{PMS} in the setting of the multi-species Sherrington-Kirkpatrick model (studied previously e.g. in \cite{FPV0}, \cite{MS}) will yield the above behaviour of the overlaps.

Notice that for $\sigma\in\Sigma(d)$, the last row and column of the overlap matrix $R_{\ell,\ell'}$ are determined by the $(\kappa-1)\times(\kappa-1)$ principal submatrix and, for $k\leq \kappa-1,$
\begin{align}
R_{\ell,\ell'}^{k,\kappa} &= d_k - \sum_{k'\leq \kappa-1} R_{\ell,\ell'}^{k,k'},\nonumber\\
R_{\ell,\ell'}^{\kappa,k} &= d_k - \sum_{k'\leq \kappa-1} R_{\ell,\ell'}^{k',k},\label{Lift}\\
R_{\ell,\ell'}^{\kappa,\kappa} &= d_{\kappa} - \sum_{k\leq \kappa-1}d_k + \sum_{k,k'\leq \kappa-1} R_{\ell,\ell'}^{k,k'}.\nonumber
\end{align}
By symmetry, for a matrix $\gamma = (\gamma_{k,k'})_{k,k'\leq \kappa}\in\Gamma_\kappa$ these equations can be written as
\begin{align}
\gamma_{\kappa,k} = \gamma_{k,\kappa} &= d_k - \sum_{k'\leq \kappa-1} \gamma_{k,k'}\,\,
\mbox{ for } k\leq \kappa-1,
\nonumber\\
\gamma_{\kappa,\kappa} &= d_{\kappa} - \sum_{k\leq \kappa}d_k + \sum_{k,k'\leq \kappa-1} \gamma_{k,k'}.\label{LiftG}
\end{align}
As a result, we can require that functions $\pi\in \Pi_d$ take values in $\Gamma_\kappa$ subject to these constraints. We will denote such matrices by $\Gamma_\kappa(d)\subseteq \Gamma_\kappa.$

Functionals of $\pi\in \Pi_d$ that will appear in the computation of the free energy will be Lipschitz with respect to the metric 
\begin{equation}
\Delta(\pi,\pi') = \int_{0}^1 \bigl\|\pi(x)-\pi'(x)\bigr\|_1 \,dx
\label{metric1}
\end{equation}
where $\|\gamma\|_1 = \sum_{k,k'} |\gamma_{k,k'}|,$ and we will explain in Section \ref{Sec4label} that a general $\pi\in\Pi_d$ can be easily discretized in a way that approximates $\pi$ in this metric. For some $r\geq 1$, a discrete path in $\Pi_d$ can be encoded by two sequences
\begin{equation}
x_{-1}=0\leq x_0 \leq \ldots \leq x_{r-1} \leq x_r = 1
\label{xs}
\end{equation}
and a monotone sequence of Gram matrices in $\Gamma_\kappa(d)$,
\begin{equation}
0=\gamma_0\leq \gamma_1 \leq \ldots \leq \gamma_{r-1}\leq \gamma_r = \diag(d_1,\ldots,d_\kappa).
\label{gammas}
\end{equation}
We can associate with these sequences a path defined by 
\begin{equation}
\pi(x)= \gamma_p \,\,\mbox{ for }\,\, x_{p-1}<x\leq x_p
\label{gammapath}
\end{equation}
for $0\leq p\leq r$, with $\pi(0)=0$. Given such a discrete path, let us consider a sequence of independent Gaussian vectors $z_p = (z_p(k))_{k\leq \kappa}$ for $0\leq p\leq r$ with the covariances
\begin{equation}
\Cov(z_p) = 2(\gamma_p - \gamma_{p-1}).
\label{zpseq}
\end{equation}
Given $\lambda = (\lambda_k)_{k\leq \kappa-1} \in \Reals^{\kappa-1}$ (which will play the role of Lagrange multipliers for the constraints in (\ref{SigmaD})), let us define 
\begin{equation}
X_r
=
 \log \sum_{k\leq \kappa} \exp\Bigl(
\beta \sum_{1\leq p \leq r} z_p(k) +\lambda_{k} \I(k\leq \kappa-1)\Bigr)
\label{Xr}
\end{equation}
(keeping the dependence of $X_r$ on $\lambda$ implicit) and, recursively over $0\leq p\leq r-1,$ define
\begin{equation}
X_p=\frac{1}{x_p}\log \e_p \exp x_p X_{p+1},
\label{Xp}
\end{equation}
where $\e_p$ denotes the expectation with respect to $z_{p+1}$ only. If $x_p=0$, we interpret this equation  as $X_p = \e_p X_{p+1}.$ Notice that $X_0$ is non-random, and we will denote it by
\begin{equation}
\Phi(\lambda, d, r, x, \gamma) := X_0,
\label{Phi}
\end{equation}
making the dependence on all the parameters explicit (the dependence on $d$ here is through the last constraint in (\ref{gammas})). Finally, we define the functional
\begin{align}
\PP(\lambda,d,r, x, \gamma)
= & \,\,\,
 \Phi(\lambda,d,r, x, \gamma) -\sum_{k\leq \kappa-1} \lambda_k d_k
 \nonumber
\\
&
-\frac{\beta^2}{2}
\sum_{0\leq p\leq r-1} x_p \bigl( \|\gamma_{p+1}\|_{HS}^2 - \|\gamma_{p}\|_{HS}^2\bigr),
\label{Ppar}
\end{align}
where $\|A\|_{HS}$ denotes the Hilbert-Schmidt norm of the matrix $A=(a_{ij})$, that is $\|A\|_{HS}^2 = \sum_{i,j} a_{ij}^2.$ The following is our main result.
\begin{theorem}\label{ThFE}
For any $\kappa\geq 1$, the limit of the free energy is given by
\begin{equation}
\lim_{N\to\infty} F_N
= 
\sup_{d} \inf_{\lambda,r, x,\gamma} \PP(\lambda,d,r, x, \gamma).
\label{Parisi}
\end{equation}
\end{theorem}
The formula (\ref{Parisi}) is the analogue of the classical Parisi formula \cite{Parisi79, Parisi} for the free energy in the Sherrington-Kirkpatrick model, which was first proved for general mixed even $p$-spin models in \cite{TPF}, and for general mixed $p$-spin models in \cite{PPF}. Let us make several remarks about Theorem \ref{ThFE}.

\medskip
\noindent
\emph{Remark 1.}
As in the setting of the classical Sherrington-Kirkpatrick model, one can observe that the functional (\ref{Phi}) depends on $(r,x,\gamma)$ only through the path $\pi$ in (\ref{gammapath}), so we can denote it by $\Phi(\lambda, d, \pi)$. Moreover, we will show (for more general family of functionals) that, for discrete paths, $\Phi$ is Lipschitz with respect to the metric (\ref{metric1}) and we can extend it by continuity to all $\pi\in \Pi_d$. Also, rearranging the terms, we can rewrite
\begin{align}
-\sum_{0\leq p\leq r-1} x_p \bigl( \|\gamma_{p+1}\|_{HS}^2 - \|\gamma_{p}\|_{HS}^2\bigr)
&=
- \|\gamma_{r}\|_{HS}^2 +  \sum_{1\leq p\leq r} (x_p-x_{p-1}) \|\gamma_{p}\|_{HS}^2
\nonumber
\\
&=
-\sum_{k\leq \kappa} d_{k}^2 + \int_{0}^1\! \|\pi(x)\|_{HS}^2\, dx
\label{rearrange}
\end{align}
and, therefore, we can rewrite (\ref{Ppar}) as
\begin{equation}
\PP(\lambda,d,\pi)
=
 \Phi(\lambda,d,\pi) 
-\sum_{k\leq \kappa-1} \lambda_k d_k 
- \frac{\beta^2}{2}\sum_{k\leq \kappa} d_{k}^2 
+\frac{\beta^2}{2}\int_{0}^1\! \|\pi(x)\|_{HS}^2\, dx.
\label{PPldpi}
\end{equation}
The last term can also be extended by continuity to all $\pi\in\Pi_d$ and the formula for the free energy can be rewritten in terms of $\pi$ (such order parameter appeared in the physics literature in \cite{FPV}),
\begin{equation}
\lim_{N\to\infty} F_N
= 
\sup_{d} \inf_{\lambda,\pi} \PP(\lambda,d,\pi).
\end{equation}

\medskip
\noindent
\emph{Remark 2.} Since the matrices $\gamma_p$ in (\ref{gammas}) represent possible values of the overlap matrix $R_{\ell,\ell'}$ and the overlaps in (\ref{Rllkk}) are non-negative, we can restrict $\gamma_p$ to the set $\Gamma_\kappa^+$ of Gram matrices with non-negative coefficients (also satisfying the constraints (\ref{LiftG})). We will see in the proof that the upper bound holds for more general sequences as above, without this additional structural information. However, from the proof of the lower bound it will be clear that it is enough to restrict the variational problem (\ref{Parisi}) to this subclass of matrices that `remembers' the structural properties of the original overlaps matrices.

\medskip
\noindent
\emph{Remark 3.} In the case of the classical Sherrington-Kirkpatrick model $\kappa=2$, the representation (\ref{Parisi}) differs from the usual Parisi formula, although it is equivalent. Up to a transformation, it essentially corresponds to the maximization over the free energies of subsystems with constrained magnetization $N^{-1}\sum_{i\leq N}\sigma_i.$ In the case $\kappa\geq 3$, it seems important to constrain the state sizes first, and it would be very interesting to know if one can remove $\sup_d$ in (\ref{Parisi}), for example, by showing that it is achieved on the configurations with equal group sizes.

\medskip
\noindent
\emph{Remark 4.} In addition to allowing $g_{ij}$'s to have non-zero mean of order $1/\sqrt{N}$, one can introduce some general external field term to the model, or consider a mixed $p$-spin version with $p$ spins interacting through $\I(\sigma_1=\ldots=\sigma_p).$ These modifications require only minor changes in the proof, so we do not consider them for simplicity of notation. 

\medskip
\noindent
\emph{Remark 5.}
The Potts spin glass model resembles, but is different from the Ghatak-Sherrington model \cite{GS} (which was solved rigorously in \cite{PGS} and is included as a special case in \cite{PVS}) where the spins take more than two values but interact through the product $\sigma_i\sigma_j$ as in (\ref{SKH}).

\medskip
\noindent
\emph{Remark 6.}  In \cite{TalUltra} (see also Section 15.7 in \cite{SG2}), Talagrand considered a system consisting of multiple copies of the Sherrington-Kirkpatrick model (or mixed even $p$-spin models), possibly at different temperatures, coupled by constraining the overlaps between them. He proved a natural generalization of the Guerra replica symmetry breaking upper bound \cite{Guerra} for such systems and asked whether this bound can be improved. This problem can be viewed as a special case of the vector version of the Sherrington-Kirkpatrick model mentioned above (more precisely, its mixed $p$-spin analogue), so the synchronization mechanism developed in this paper can be used to show that the bound of Talagrand is, in fact, asymptotically sharp (see \cite{PVS}). 

\medskip
We begin in Section \ref{Sec2label} with the analogue of Guerra's replica symmetry breaking interpolation and the proof of the upper bound. One of the functionals arising in this interpolation does not automatically decouple over spin coordinates, and in Section \ref{Sec3label} we prove a basic large deviations result that gives the right form of decoupling, which is used later in the proof of the lower bound. In preparation for the proof of the lower bound, in Section \ref{Sec4label} we collect various basic continuity and approximation properties of this and other functionals. Sections \ref{Sec5label} and \ref{Sec6label} contain the core new ideas of the proof. In Section \ref{Sec5label}, we prove a new family of the Ghirlanda-Guerra identities via a small perturbation of the Hamiltonian and in Section \ref{Sec6label} we utilize these identities to prove strong synchronization properties for the blocks of the overlap array that were mentioned above. Finally, we put all the pieces together in Section \ref{Sec7label}, where we prove the matching lower bound via the standard cavity computation.

\medskip
\noindent
\textbf{Acknowledgements.} The author would like to thank Giorgio Parisi for valuable comments.

\section{Upper bound via Guerra's interpolation}\label{Sec2label}

In this section we will show that the functional $\PP(\lambda,d,r, x, \gamma)$ in (\ref{Ppar}) is an upper bound for $F_N(d)$ for any $\lambda,r, x, \gamma$, using the analogue of Guerra's interpolation \cite{Guerra}.  Without loss of generality, we can and will assume that the inequalities in (\ref{xs}) are strict,
\begin{equation}
x_{-1}=0< x_0 < \ldots < x_{r-1} < x_r = 1.
\label{xsstrict}
\end{equation}
Let $(v_\alpha)_{\alpha\in \Natural^r}$ be the weights of the Ruelle probability cascades \cite{Ruelle} corresponding to the sequence (\ref{xsstrict}) (see e.g. Section 2.3 in \cite{SKmodel} for the definition). For $\alpha^1, \alpha^2\in \Natural^r$, we denote
\begin{equation}
\alpha^1\wedge \alpha^2 = \min\Bigl\{0\leq p \leq r  \mid  \alpha_1^1= \alpha_1^2, \ldots, \alpha_{p}^1 = \alpha_{p}^2, \alpha_{p+1}^1 \not = \alpha_{p+1}^2 \Bigr\},
\end{equation}
where $\alpha^1\wedge \alpha^2 =r$ if $\alpha^1=\alpha^2$. Since the sequence in (\ref{gammas}) is non-decreasing, the sequence $\|\gamma_p\|_{HS}$ is also non-decreasing. As a result, there exist Gaussian processes 
\begin{equation}
Z^\alpha=(Z^\alpha(k))_{k\leq\kappa} \mbox{ and } Y^\alpha,
\end{equation}
both indexed by $\alpha\in\Natural^r$, with the covariances
\begin{align}
\Cov(Z^{\alpha^1}, Z^{\alpha^2}) &= 2\gamma_{\alpha^1\wedge\alpha^2},
\nonumber
\\
\Cov(Y^{\alpha^1}, Y^{\alpha^2}) &= \|\gamma_{\alpha^1\wedge\alpha^2}\|_{HS}^2.
\label{CD}
\end{align}
Let $Z_i^\alpha$ be independent copies of the process $Z^\alpha$, also independent of $Y^\alpha$. For $0\leq t\leq 1$, consider an interpolating Hamiltonian defined on $\Sigma_N\times \Natural^r$ by
\begin{equation}
H_{N,t}(\sigma,\alpha) = 
\sqrt{t} H_N(\sigma) + \sqrt{1-t} \sum_{i\leq N} Z_{i}^\alpha(\sigma_i) +\sqrt{t} \sqrt{N} Y^\alpha.
\label{HNt}
\end{equation}
Similarly to (\ref{FED}), we define
\begin{equation}
\varphi(t):= \frac{1}{N} \e \log \sum_{\alpha\in\Natural^r} v_\alpha \sum_{\sigma\in \Sigma(d)} \exp \beta H_{N,t}(\sigma,\alpha).
\label{FEDt}
\end{equation}
Then it is easy to check the following.
\begin{lemma}\label{LemGUP}
The function $\varphi(t)$ in (\ref{FEDt}) is non-increasing.
\end{lemma}
\textbf{Proof.} Let us denote by $\la\, \cdot\,\ra_t$ the average with respect to the Gibbs measure 
$$
G_t(\sigma,\alpha) \sim v_\alpha \exp \beta H_{N, t}(\sigma,\alpha).
$$
on $\Sigma(d)\times \Natural^r$. Then, for $0<t<1$,
$$
\varphi'(t) = \frac{\beta}{N}\e \Bigl\la \frac{\partial H_{N, t}(\sigma,\alpha)}{\partial t} \Bigr\ra_t.
$$
If we rewrite the Hamiltonian $H_N(\sigma)$ as
\begin{equation}
H_N(\sigma) = \frac{1}{\sqrt{N}} \sum_{k=1}^{\kappa}\sum_{i,j=1}^N g_{ij}\I(\sigma_i = k)\I(\sigma_j = k)
\label{DeanHk}
\end{equation}
and recall the definition of the overlaps in (\ref{Rll}), a direct calculation gives
\begin{equation}
\Cov\bigl(H_N(\sigma^1), H_N(\sigma^2)\bigr) = N\sum_{k,k'\leq \kappa}(R_{1,2}^{k,k'})^2.
\label{CovHNsum}
\end{equation}
Similarly, if we write $Z_i^\alpha(\sigma_i) = \sum_{k\leq \kappa}\I(\sigma_i=k)Z_i^\alpha(k)$ then, from the definition (\ref{CD}),
\begin{equation}
\Cov\Bigl(\sum_{i\leq N} Z_{i}^{\alpha^1}(\sigma_i^1), \sum_{i\leq N} Z_{i}^{\alpha^2}(\sigma_i^2)\Bigr)
=
2N\sum_{k,k'\leq \kappa} R_{1,2}^{k,k'} \gamma_{\alpha^1\wedge\alpha^2}^{k,k'},
\end{equation}
where $\gamma_{\alpha^1\wedge\alpha^2}^{k,k'}$ is the $(k,k')$-element of the matrix $ \gamma_{\alpha^1\wedge\alpha^2}.$ Using these equations and recalling the covariance of $Y^\alpha$ in (\ref{CD}),
$$
\frac{1}{N}\, \e \frac{\partial H_{N, t}(\sigma^1,\alpha^1)}{\partial t} H_{N, t}(\sigma^2,\alpha^2)
=
\frac{1}{2}\sum_{k,k'\leq \kappa}
\bigl(
R_{1,2}^{k,k'} - \gamma_{\alpha^1\wedge\alpha^2}^{k,k'} 
\bigr)^2.
$$
For $(\sigma^1,\alpha^1)=(\sigma^2,\alpha^2)$, this sum vanishes because $R_{1,1} = \gamma_r = \diag(d_1,\ldots,d_{\kappa})$. Finally, usual Gaussian integration by parts (see e.g. Lemma 1.1 in \cite{SKmodel}) yields
$$
\varphi'(t) = - \frac{\beta^2}{2}\sum_{k,k'\leq \kappa}  
\e \Bigl\la \bigl(
R_{1,2}^{k,k'} - \gamma_{\alpha^1\wedge\alpha^2}^{k,k'} 
\bigr)^2
\Bigr\ra_t \leq 0.
$$
This finishes the proof.
\qed

\medskip
The lemma implies that $\varphi(1)\leq \varphi(0).$ First of all,
$$
\varphi(1)
=
F_N(d)
+
\frac{1}{N}\smsp \e\log \sum_{\alpha\in\Natural^r} v_{\alpha} 
\exp \beta \sqrt{N} Y^\alpha.
$$
The standard properties of the Ruelle probability cascades (see Section 2.3 and the proof of Lemma 3.1 in \cite{SKmodel}) together with the covariance structure (\ref{CD}) imply that 
\begin{equation}
\frac{1}{N}\smsp \e\log \sum_{\alpha\in\Natural^r} v_{\alpha} 
\exp \beta \sqrt{N} Y^\alpha
=
\frac{\beta^2}{2}
\sum_{0\leq p\leq r-1} x_p \bigl( \|\gamma_{p+1}\|_{HS}^2 - \|\gamma_{p}\|_{HS}^2\bigr).
\label{funp2}
\end{equation}
Next, let us consider
$$
\varphi(0)
=
\frac{1}{N} \e \log \sum_{\alpha\in\Natural^r} v_\alpha \sum_{\sigma\in \Sigma(d)} \exp \beta \sum_{i\leq N} Z_{i}^\alpha(\sigma_i).
$$
Since $\sum_{i\leq N} \I(\sigma_i = k) = Nd_k$ for $\sigma\in\Sigma(d)$, adding $\sum_{k\leq \kappa-1} \lambda_k \sum_{i\leq N} \I(\sigma_i = k)$ and at the same time subtracting $N\sum_{k\leq \kappa-1} \lambda_k d_k$ for any $\lambda_k\in\Reals$ in the exponent will not change $\varphi(0)$. If we then replace the sum over $\sigma\in \Sigma(d)$ by the sum over all $\sigma$, we obtain the upper bound
$$
\varphi(0)
\leq
-\sum_{k\leq \kappa-1} \lambda_k d_k +
\frac{1}{N} \e \log \sum_{\alpha\in\Natural^r} v_\alpha \sum_{\sigma} \exp \sum_{i\leq N}
\Bigl( \beta Z_{i}^\alpha(\sigma_i) + \sum_{k\leq \kappa-1} \lambda_k \I(\sigma_i = k)\Bigr).
$$
If we introduce the notation
$$
X_i^\alpha =
\sum_{\sigma_i\leq \kappa}\exp\Bigl( \beta Z_{i}^\alpha(\sigma_i) 
+ \sum_{k\leq \kappa-1} \lambda_k \I(\sigma_i = k)\Bigr)
=
\sum_{\sigma_i\leq \kappa}\exp\Bigl( \beta Z_{i}^\alpha(\sigma_i) 
+ \lambda_{\sigma_i} \I(\sigma_i \leq \kappa-1)\Bigr)
$$
then this upper bound can be rewritten as
$$
\varphi(0)
\leq
-\sum_{k\leq \kappa-1} \lambda_k d_k +
\frac{1}{N} \e \log \sum_{\alpha\in\Natural^r} v_\alpha 
\prod_{i\leq N} X_{i}^\alpha.
$$
Again, standard properties of the Ruelle probability cascades (see Section 2.3 in \cite{SKmodel}) imply that
$$
\frac{1}{N} \e \log \sum_{\alpha\in\Natural^r} v_\alpha \prod_{i\leq N} X_{i}^\alpha
=
\e \log \sum_{\alpha\in\Natural^r} v_\alpha X_{1}^\alpha = X_0,
$$
where $X_0 = \Phi(\lambda, d, r, x, \gamma)$ was defined in (\ref{Phi}) and, therefore,
\begin{equation}
\varphi(0) \leq -\sum_{k\leq \kappa-1} \lambda_k d_k + \Phi(\lambda, d, r, x, \gamma).
\label{phi0bound}
\end{equation}
Together with the inequality $\varphi(1)\leq \varphi(0)$ this implies that $F_N(d) \leq \PP(\lambda,d,r, x, \gamma)$, which proves the upper bound in Theorem \ref{ThFE}.

\section{Decoupling the constraints on sizes of states}\label{Sec3label}

The quantity $\varphi(0)$ will also appear in the proof of the lower bound and, at that moment, we will need to use the fact that the upper bound (\ref{phi0bound}) becomes asymptotically exact after we minimize over $\lambda = (\lambda_k)_{k\leq \kappa-1},$ which we will prove in this section. This type of feature first appeared in the spherical Sherrington-Kirkpatrick model (see \cite{T-sphere, Chen-sphere}), as well as in the study of the Ghatak-Sherrington model in \cite{PGS}.

In the notation of the previous section, let 
\begin{equation}
f_N(d):= \frac{1}{N} \e \log \sum_{\alpha\in\Natural^r} v_\alpha \sum_{\sigma\in \Sigma(d)} \exp \beta \sum_{i\leq N} Z_{i}^\alpha(\sigma_i).
\label{fN}
\end{equation}
Let us also note right away that in the computation leading to (\ref{phi0bound}) we showed that
\begin{equation}
\Phi(\lambda, d, r, x, \gamma) 
=
\frac{1}{N} \e \log \sum_{\alpha\in\Natural^r} v_\alpha \sum_{\sigma} \exp \sum_{i\leq N}
\Bigl( \beta Z_{i}^\alpha(\sigma_i) + \sum_{k\leq \kappa-1} \lambda_k \I(\sigma_i = k)\Bigr)
\label{PhiLrer}
\end{equation}
for any $N$. Let us consider the set
\begin{equation}
\DD_N = \Bigl\{d\in\DD \mid \Sigma(d) \mbox{ is not empty}\Bigr\}.
\label{DDN}
\end{equation}
We will now prove the following.
\begin{lemma}\label{Lem2}
If $d^N\in \DD_N$ and $\lim_{N\to\infty} d^N = d$ then
\begin{equation}
\lim_{N\to\infty} f_N(d^N) = \inf_{\lambda}\Bigl(
-\sum_{k\leq \kappa-1} \lambda_k d_k + \Phi(\lambda, d, r, x, \gamma)
\Bigr).
\label{fNlim}
\end{equation}
\end{lemma}
Before we begin the proof, let us point out one subtle point. Notice that $f_N(d)$ depends on $d$ through the constraint $\sigma\in\Sigma(d)$, but also through the covariance structure of $Z_i^\alpha$ due to the last constraint in (\ref{gammas}). The dependence of $\Phi$ on $d$ is only through this covariance structure. For the rest of this section, we will fix this covariance structure so that the dependence on $d$ will be only through the constraint $\sigma\in\Sigma(d)$. In other words, we will prove (\ref{fNlim}) even if the covariance structure is not related to the constraint on $\sigma.$ In particular, since $r,x,\gamma$ are also fixed, we will view $\Phi = \Phi(\lambda)$ as a function of $\lambda$ only and show that
\begin{equation}
\lim_{N\to\infty} f_N(d^N) = \inf_{\lambda}\Bigl(
-\sum_{k\leq \kappa-1} \lambda_k d_k + \Phi(\lambda)
\Bigr).
\label{fNlim2}
\end{equation}
\textbf{Proof of Lemma \ref{Lem2}.} We will first give an outline of the proof, and all the steps will be completed in the rest of the section. In the first step we will show that, for all $d\in \DD$, the limit 
\begin{equation}
f(d) = \lim_{N\to\infty} f_N(d^N)
\label{flim}
\end{equation}
exists and is concave in $d$. Since, the function $f_N(d)$ is, clearly, bounded from above and below uniformly over $N$ and $d$ such that $\Sigma(d)$ is not empty, the limit $f(d)$ will be bounded and continuous on $\DD$. In the second step we will show that 
\begin{equation}
\Phi(\lambda) = \max_{d\in\DD}\Bigl(f(d) + \sum_{k\leq \kappa-1}\lambda_k d_k\Bigr).
\label{LFmax}
\end{equation}
By the biconjugation theorem for convex functions (see e.g. Theorem 12.2 in \cite{Rockafellar}), we will then conclude that (\ref{fNlim2}) holds.
\qed

\medskip
In order to prove (\ref{flim}), instead of working with the sequence $(d^N)$ it will be convenient to relax the constraint $\sigma\in\Sigma(d)$ instead. Given $\eps>0$, we define
\begin{equation}
\Sigma_\eps(d) =  \Bigl\{ \sigma \mid \sum_{i\leq N}\I(\sigma_i=k)\in N[d_k-\eps,d_k+\eps] \,\mbox{ for all }\, k\leq \kappa \Bigr\}
\label{SigmaDE}
\end{equation}
and, similarly to (\ref{fN}), define
\begin{equation}
f_{N,\eps}(d):= \frac{1}{N} \e \log \sum_{\alpha\in\Natural^r} v_\alpha \sum_{\sigma\in \Sigma_\eps(d)} \exp \beta \sum_{i\leq N} Z_{i}^\alpha(\sigma_i).
\end{equation}
We begin with the following observation (which is an adaptation of Lemma 1 in \cite{SKcoupled}).
\begin{lemma}\label{Lem3}
There exists a constant $L>0$ independent of $N$ such that
\begin{equation}
\sup_{d\in\DD_N}|f_{N,\eps}(d)-f_N(d)|\leq L\sqrt{\eps}.
\label{Lem3eq}
\end{equation}
\end{lemma}
In particular, since, for any $d^1,d^2\in\DD_N$ and $\eps=\max_{k\leq \kappa}|d_k^1-d_k^2|$, we have the inclusions
$$
\Sigma(d^1)\subseteq \Sigma_\eps(d^2),\,\, \Sigma(d^2)\subseteq \Sigma_\eps(d^1),
$$ 
Lemma \ref{Lem3} implies that
\begin{equation}
|f_N(d^1) - f_N(d^2)| \leq L\max_{k\leq \kappa}|d_k^1-d_k^2|^{1/2} \leq L\|d^1-d^2\|_\infty^{1/2}.
\label{dinft}
\end{equation}

\medskip
\noindent
\textbf{Proof of Lemma \ref{Lem3}.} For $\sigma\in\Sigma_\eps(d)$, let $\tilde{\sigma}$ be a vector in $\Sigma(d)$ with the smallest number of different coordinates $\sum_{i\leq N}\I(\sigma_i\not= \tilde{\sigma}_i)$. Then it is obvious that
$
\sum_{i\leq N}\I(\sigma_i\not= \tilde{\sigma}_i) \leq L N\eps
$
for some constant $L$ that depends on $\kappa$ only. First, we will compare $f_{N,\eps}(d)$ with 
$$
\tilde{f}_{N,\eps}(d):= \frac{1}{N} \e \log \sum_{\alpha\in\Natural^r} v_\alpha \sum_{\sigma\in \Sigma_\eps(d)} \exp \beta \sum_{i\leq N} Z_{i}^\alpha(\tilde{\sigma}_i).
$$
Let $\tilde{Z}_i^\alpha$ be independent copies of the processes $Z_i^\alpha$, let
$$
Z_t(\alpha,\sigma) = \sum_{i\leq N} \bigl( \sqrt{t }Z_{i}^\alpha(\sigma_i)+ \sqrt{1-t}\tilde{Z}_{i}^\alpha(\tilde{\sigma}_i)\bigr)
$$
for $t\in [0,1]$, and consider the interpolation
$$
\varphi(t)= \frac{1}{N} \e \log \sum_{\alpha\in\Natural^r} v_\alpha \sum_{\sigma\in \Sigma_\eps(d)} \exp \beta Z_t(\alpha,\sigma)
$$
such that $\varphi(1) = {f}_{N,\eps}(d)$ and $\varphi(0)=\tilde{f}_{N,\eps}(d)$. One can compute the derivative $\varphi'(t)$ using Gaussian integration by parts as in Lemma \ref{LemGUP}. If we recall the covariance formulas in (\ref{CD}),
$$
\frac{1}{N} \e \frac{\partial Z_t(\alpha^1,\sigma^1)}{\partial t} Z_{t}(\alpha^2,\sigma^2)
=
\frac{1}{N}  \sum_{i\leq N} \Bigl( \gamma_{\alpha^1\wedge\alpha^2}^{\sigma_i^1,\sigma_i^2}-\gamma_{\alpha^1\wedge\alpha^2}^{\tilde{\sigma}_i^1,\tilde{\sigma}_i^2} \Bigr).
$$
The $i^{\mathrm{th}}$ term is zero unless $\sigma_i^1\not = \tilde{\sigma}_i^1$ or $\sigma_i^2\not = \tilde{\sigma}_i^2$ and, by the definition of $\tilde{\sigma}$ above, the number of such coordinates is bounded by $LN\eps.$ Therefore, $|\varphi'(t)|\leq L\beta^2\eps$ and 
$$
|{f}_{N,\eps}(d)- \tilde{f}_{N,\eps}(d)| \leq L\beta^2\eps
$$
for some constant $L$ that depends only on $\kappa$.

For $\sigma\in \Sigma(d)$, let us denote by $\N(\sigma)$ the number of configurations $\rho\in\Sigma_\eps(d)$ such that $\tilde{\rho} = \sigma.$ Then we can rewrite and bound $\tilde{f}_{N,\eps}(d)$ as follows,
\begin{align*}
\tilde{f}_{N,\eps}(d) 
&= 
\frac{1}{N} \e \log \sum_{\alpha\in\Natural^r} v_\alpha \sum_{\sigma\in \Sigma(d)} \N(\sigma)\exp \beta \sum_{i\leq N} Z_{i}^\alpha({\sigma}_i).
\\
&\leq f_{N}(d) +\frac{1}{N}\max_{\sigma\in \Sigma(d)} \log \N(\sigma).
\end{align*}
For any $\sigma\in \Sigma(d)$, the number $\N(\sigma)$ is bounded by the number of configurations $\rho$ such that $\sum_{i\leq N}\I(\rho_i\not= {\sigma}_i) \leq L N\eps$. By the classical large deviation estimate for Bernoulli random variables, a number of different ways to choose $LN\eps$ coordinates is bounded by  $2^N \exp(-N I(1-L\eps)),$ where 
$$
I(x)=\frac{1}{2}\bigl((1+x)\log(1+x)+(1-x)\log(1-x)\bigr),
$$
and there are $\kappa^{LN\eps}$ ways to choose $\rho_i$'s different from $\sigma_i$'s on these coordinates. Therefore,
\begin{align*}
\frac{1}{N}\max_{\sigma\in \Sigma(d)} \log \N(\sigma)
&\leq
L\eps\log\kappa+
\log 2 - I(1-L\eps)
\\
&=
L\eps\log\kappa+
\log\Bigl(1+\frac{L\eps}{2-L\eps}\Bigr)+\frac{L\eps}{2}\log
\frac{2-L\eps}{\eps}
\leq L\sqrt{\eps},
\end{align*}
for small enough $\eps$. We showed that ${f}_{N,\eps}(d)\leq f_N(d)+L\sqrt{\eps}$ and, since $ f_N(d)\leq {f}_{N,\eps}(d),$ this finishes the proof.
\qed

\medskip
\begin{lemma}
For any $d\in\DD$, the limit
\begin{equation}
\lim_{N\to\infty} f_{N,\eps}(d) = f_\eps(d)
\label{fepsd}
\end{equation}
exists and is a concave function of $d$.
\end{lemma}
\textbf{Proof.}
Let us make the dependence of the set $\Sigma^N_\eps(d)$ in (\ref{SigmaDE}) on $N$ explicit. For any $N_1,N_2\geq 1$, let $N=N_1+N_2$ and $\lambda = N_1/N.$ For any $d^1,d^2\in\DD$, let $d=\lambda d^1 + (1-\lambda)d^2.$ Then, clearly,
$$
\Sigma^N_\eps(d)\supseteq \Sigma^{N_1}_\eps(d^1)\times \Sigma^{N_2}_\eps(d^2)
$$
and, therefore, the following inequality holds,
$$
N f_{N,\eps}(d)\geq N_1 f_{N_1,\eps}(d^1)+N_2 f_{N_2,\eps}(d^2). 
$$
In particular, for $d^1=d^2=d$ this shows that the sequence $N f_{N,\eps}(d)$ is super-additive and, hence, the limit (\ref{fepsd}) exists. Dividing both sides by $N$, we get
$$
f_{N,\eps}(d)= f_{N,\eps}(d^1 + (1-\lambda)d^2)\geq \lambda f_{N_1,\eps}(d^1)+(1-\lambda) f_{N_2,\eps}(d^2)
$$
and taking limits shows that $f_\eps(d)$ is concave.
\qed

\medskip
\noindent
Combining the above two lemmas, we complete the first step. 
\begin{lemma}\label{Lem5}
If $d^N\in \DD_N$ and $\lim_{N\to\infty} d^N = d\in \DD$ then the limit
\begin{equation}
f(d):= \lim_{N\to\infty} f_N(d^N) = \lim_{\eps\downarrow 0} f_\eps(d)
\label{Lem5eq1}
\end{equation}
exists and, for all $d^1,d^2\in\DD$, satisfies
\begin{equation}
|f(d^1) - f(d^2)| \leq L\|d^1-d^2\|_\infty^{1/2}.
\label{Lem5eq2}
\end{equation}
\end{lemma}
\textbf{Proof.}
Suppose that $\delta := \max_{k\leq \kappa}|d_k^N-d_k|\leq \eps.$ The inclusions
$\Sigma(d^N)\subseteq \Sigma_\eps(d)\subseteq \Sigma_{\eps+\delta}(d^N)$
together with (\ref{Lem3eq}) imply that
$$
f_N(d^N) \leq f_{N,\eps}(d)\leq f_{N,\eps+\delta}(d^N) \leq f_N(d^N) + L\sqrt{\eps+\delta}.
$$
Taking limits and using (\ref{fepsd}) implies that
$$
\limsup_{N\to\infty} f_N(d^N) \leq f_\eps(d)\leq \liminf_{N\to\infty} f_N(d^N) + L\sqrt{\eps}.
$$
Finally, letting $\eps\downarrow 0$ proves (\ref{Lem5eq1}). By (\ref{dinft}), this limit satisfies (\ref{Lem5eq2}).
\qed

\medskip
Next, we will focus on the second step (\ref{LFmax}). 
\begin{lemma} For any $\lambda = (\lambda_k)_{k\leq\kappa-1}\in\Reals^{\kappa-1}$,
\begin{equation}
\Phi(\lambda) = \max_{d\in\DD}\Bigl(f(d) + \sum_{k\leq \kappa-1}\lambda_k d_k\Bigr).
\end{equation}
\end{lemma}
\textbf{Proof.}
This will be a direct consequence of the properties of the Ruelle probability cascades. More precisely, we will use a recursive representation for  functionals of the type
\begin{equation}
\Phi(\lambda, S) 
=
\frac{1}{N} \e \log \sum_{\alpha\in\Natural^r} v_\alpha \sum_{\sigma\in S} \exp \sum_{i\leq N}
\Bigl( \beta Z_{i}^\alpha(\sigma_i) + \sum_{k\leq \kappa-1} \lambda_k \I(\sigma_i = k)\Bigr)
\label{PhiLS}
\end{equation}
for non-empty subsets $S\subseteq \Sigma = \{1,\ldots,\kappa\}^N$. Notice, for example, that, by (\ref{fN}) and (\ref{PhiLrer}),
$$
\Phi(\lambda, \Sigma) = \Phi(\lambda)
\,\mbox{ and }\,
\Phi(\lambda, \Sigma(d)) = f_N(d) + \sum_{k\leq \kappa-1} \lambda_k d_k.
$$ 
Let us recall the definition of the sequence $z_p = (z_p(k))_{k\leq \kappa}$ for $0\leq p\leq r$ in (\ref{zpseq}) and let $(z_{p,i})_{0\leq p\leq r}$ be its independent copies for $i\leq N.$ Let us define 
\begin{equation}
X_r(\lambda, S)
=
 \log \sum_{\sigma\in S} 
 \exp \sum_{i\leq N} \Bigl(
\beta \sum_{1\leq p \leq r} z_{p,i}(\sigma_i) +\sum_{k\leq \kappa-1} \lambda_{k} \I(\sigma_i=k)\Bigr)
\label{XrLS}
\end{equation}
and, recursively over $0\leq p\leq r-1,$ define
\begin{equation}
X_p(\lambda, S)=\frac{1}{x_p}\log \e_p \exp x_p X_{p+1}(\lambda, S),
\label{XpLS}
\end{equation} 
where $\e_p$ denotes the expectation with respect to $z_{p+1,i}$ for $i\leq N$. Standard properties of the Ruelle probability cascades (Theorem 2.9 in \cite{SKmodel}) imply that
\begin{equation}
\Phi(\lambda, S) = \frac{1}{N}X_0(\lambda, S).
\end{equation}
Now notice that, since $\Sigma$ is a disjoint union of $\Sigma(d)$ for $d\in \DD_N$, 
$$
\exp X_r(\lambda, \Sigma) = \sum_{d\in\DD_N} \exp X_r(\lambda, \Sigma(d)).
$$
Since $x_{r-1}\leq 1$,
\begin{align*}
\exp x_{r-1} X_{r-1}(\lambda, \Sigma)
&=
\e_{r-1} \exp x_{r-1} X_{r}(\lambda, \Sigma)
=
\e_{r-1}\Bigl(\sum_{d\in\DD_N} \exp X_r(\lambda, \Sigma(d))\Bigr)^{x_{r-1}}
\\
&\leq 
\sum_{d\in\DD_N} \e_{r-1} \exp x_{r-1} X_r(\lambda, \Sigma(d))
=
\sum_{d\in\DD_N} \exp x_{r-1} X_{r-1}(\lambda, \Sigma(d)).
\end{align*}
By induction, using that $x_p/ x_{p+1} \leq 1,$
\begin{align*}
\exp x_{p} X_{p}(\lambda, \Sigma)
&=
\e_{p} \exp x_{p} X_{p+1}(\lambda, \Sigma)
\leq 
\e_p \Bigl( \sum_{d\in\DD_N} \exp x_{p+1} X_{p+1}(\lambda, \Sigma(d))\Bigr)^{x_p/x_{p+1}}
\\
&\leq
\sum_{d\in\DD_N} \e_{p} \exp x_{p} X_{p+1}(\lambda, \Sigma(d))
=
\sum_{d\in\DD_N} \exp x_{p} X_{p}(\lambda, \Sigma(d)).
\end{align*}
Recall that we assumed in (\ref{xsstrict}) that $x_0>0$ so, for $p=0$, this gives
\begin{align*}
\Phi(\lambda) = \frac{1}{N} X_{0}(\lambda, \Sigma)
&\leq
\frac{1}{Nx_0}
\log \sum_{d\in\DD_N} \exp x_{0} X_{0}(\lambda, \Sigma(d))
\\
&\leq
\frac{1}{Nx_0} \log \mathrm{card}(\DD_N)
+
\max_{d\in\DD_N} \frac{1}{N} X_{0}(\lambda, \Sigma(d))
\\
&\leq
\frac{1}{Nx_0} \log N^{\kappa}
+
\max_{d\in\DD_N} \Bigl(f_N(d) + \sum_{k\leq \kappa-1} \lambda_k d_k \Bigr).
\end{align*}
Combining this with (\ref{phi0bound}), we get
$$
0\leq \Phi(\lambda) -
\max_{d\in\DD_N} \Bigl(f_N(d) + \sum_{k\leq \kappa-1} \lambda_k d_k \Bigr)
\leq \frac{1}{Nx_0} \log N^{\kappa}.
$$
Using Lemma \ref{Lem5} finishes the proof.
\qed

\section{Continuity and discretization}\label{Sec4label}

In this section we will collect several technical continuity and approximation properties for various functionals that already appeared above and will appear below in the proof of the lower bound. The first two Lipschitz continuity properties are direct analogues of a well-known result of Guerra in \cite{Guerra} (see also \cite{PM} or Theorem 14.11.2 in \cite{SG2}) in the setting of the SK model.

\medskip
\noindent
\textbf{Lipschitz continuity I.} First, let us consider an arbitrary non-empty subset $S\subseteq \{1,\ldots,\kappa\}^N$ and consider the functional 
\begin{equation}
f^1_N(S,\pi):= \frac{1}{N} \e \log \sum_{\alpha\in\Natural^r} v_\alpha \sum_{\sigma\in S} \exp \beta \sum_{i\leq N} Z_{i}^\alpha(\sigma_i).
\label{fNS}
\end{equation}
defined similarly to (\ref{fN}). Right now we view this as a functional of the path $\pi$ defined in (\ref{gammapath}) in terms of the sequences (\ref{gammas}) and (\ref{xsstrict}), which determines the covariance structure of the Gaussian processes $Z_{i}^\alpha$, and we are interested in the continuity properties of $f^1_N$.

\begin{lemma}\label{LemLC1}
For two discrete paths $\pi, \tilde{\pi}\in \Pi_d$ as in (\ref{gammapath}),
\begin{equation}
\bigl|f^1_N(S,\pi)-f^1_N(S,\tilde{\pi}) \bigr|
\leq
L \int_{0}^{1}\! \|\pi(x) - \tilde{\pi}(x)\|_1 \, dx,
\label{Lipschitz1}
\end{equation}
where the constant $L$ does not depend on $N$ or the set $S$.
\end{lemma}
\textbf{Proof.} Without loss of generality, we can suppose that $\pi$, $\tilde{\pi}$ are defined in terms of the sequences
\begin{align*}
x_{-1} =0 &< x_0 < \ldots < x_{r-1} < x_r = 1,\\
0 &=\gamma_0\leq \ldots \leq \gamma_{r-1}\leq \gamma_r = \diag(d_1,\ldots,d_\kappa),\\
0 &=\tilde{\gamma}_0\leq \ldots \leq \tilde{\gamma}_{r-1}\leq \tilde{\gamma}_r = \diag(d_1,\ldots,d_\kappa).
\end{align*}
In other words, the same $x_p$s are used to define $\pi$, $\tilde{\pi}$. The reason for this is because we can combine the two sequences $(x_p)$ by artificially inserting additional values $(\gamma_p)$ and $(\tilde{\gamma}_p)$, because, as we mentioned above, the functional depends on these sequences only through $\pi$ and $\tilde{\pi}$. 

For $0\leq t\leq 1,$ the sequence $\gamma_p^t = t\gamma_p+(1-t)\tilde{\gamma}_p$ is, again, non-decreasing in $\Gamma_\kappa$. If we consider independent Gaussian processes 
$$
Z^\alpha=(Z^\alpha(k))_{k\leq\kappa},\,\, \tilde{Z}^\alpha=(\tilde{Z}^\alpha(k))_{k\leq\kappa}
$$ 
indexed by $\alpha\in\Natural^r$, with the covariances
\begin{align*}
\Cov(Z^{\alpha^1}, Z^{\alpha^2}) &= 2\gamma_{\alpha^1\wedge\alpha^2},
\\
\Cov(\tilde{Z}^{\alpha^1}, \tilde{Z}^{\alpha^2}) &= 2\tilde{\gamma}_{\alpha^1\wedge\alpha^2},
\end{align*}
then $Z^\alpha_t = \sqrt{t}Z^\alpha+\sqrt{1-t} \tilde{Z}^\alpha$ will have the covariance $\Cov(Z_t^{\alpha^1}, Z_t^{\alpha^2}) = 2\gamma^t_{\alpha^1\wedge\alpha^2}.$ Let us consider independent copies $Z^\alpha_{t,i}$ of this process for $i\leq N$ and define
\begin{equation}
\varphi(t):= \frac{1}{N} \e \log \sum_{\alpha\in\Natural^r} v_\alpha \sum_{\sigma\in S} \exp \beta \sum_{i\leq N} Z_{t,i}^\alpha(\sigma_i).
\label{fNS2}
\end{equation}
Then $\varphi(1) = f_N^1(S,\pi)$ and $\varphi(0) = f_N^1(S,\tilde{\pi}).$ To finish the proof, we will compute the derivative as in Lemma \ref{LemGUP}. Let us denote by $\la\, \cdot\,\ra_t$ the average with respect to the Gibbs measure 
\begin{equation}
G_t(\sigma,\alpha) \sim v_\alpha \exp \beta H_{N, t}(\sigma,\alpha).
\label{GtsA}
\end{equation}
on $S\times \Natural^r$ and let now $H_{N, t}(\sigma,\alpha) =  \sum_{i\leq N} Z_{t,i}^\alpha(\sigma_i).$ Then, for $0<t<1,$
$$
\varphi'(t) = \frac{\beta}{N}\e \Bigl\la \frac{\partial H_{N, t}(\sigma,\alpha)}{\partial t} \Bigr\ra_t.
$$
Recalling the covariance of $Z^\alpha$ and $\tilde{Z}^\alpha$ above, it is easy to see that
$$
\frac{1}{N}\, \e \frac{\partial H_{N, t}(\sigma^1,\alpha^1)}{\partial t} H_{N, t}(\sigma^2,\alpha^2)
=
\frac{1}{N}\sum_{i\leq N} \Bigl(
\gamma_{\alpha^1\wedge\alpha^2}^{\sigma_i^1,\sigma_i^2}
-\tilde{\gamma}_{\alpha^1\wedge\alpha^2}^{\sigma_i^1,\sigma_i^2}
\Bigr),
$$
which is zero for $(\sigma^1,\alpha^1)=(\sigma^2,\alpha^2)$ and can be bounded in absolute value by
$
\bigl \| \gamma_{\alpha^1\wedge\alpha^2}
-\tilde{\gamma}_{\alpha^1\wedge\alpha^2} \bigr\|_1,
$
where $\|\gamma\|_1 = \sum_{k,k'} |\gamma_{k,k'}|$. Therefore, Gaussian integration by parts gives
$$
|\varphi'(t)| \leq \beta^2 \e \Bigl\la \bigl \| \gamma_{\alpha^1\wedge\alpha^2}
-\tilde{\gamma}_{\alpha^1\wedge\alpha^2} \bigr\|_1 \Bigr\ra_t.
$$
For any $0\leq t\leq 1,$ the marginal of the random measure (\ref{GtsA}) on $\Natural^r$ has the same distribution as the weights $(v_\alpha)_{\alpha\in\Natural^r}$ (see e.g. Theorem 4.4 in \cite{SKmodel}) and, as a result,
\begin{align*}
 \e \bigl\la \bigl \| \gamma_{\alpha^1\wedge\alpha^2}
-\tilde{\gamma}_{\alpha^1\wedge\alpha^2} \bigr\|_1 \bigr\ra_t
&=
\e \sum_{\alpha^1,\alpha^2} v_{\alpha^1}v_{\alpha^2} \bigl \| \gamma_{\alpha^1\wedge\alpha^2}
-\tilde{\gamma}_{\alpha^1\wedge\alpha^2} \bigr\|_1
\\
&=
\sum_{0\leq p\leq r} \bigl \| \gamma_{p} -\tilde{\gamma}_{p} \bigr\|_1 \e \sum_{\alpha^1\wedge\alpha^2 = p} v_{\alpha^1}v_{\alpha^2}
\\
\mbox{(see eq. (2.82) in \cite{SKmodel}) }
&=
\sum_{0\leq p\leq r} \bigl \| \gamma_{p} -\tilde{\gamma}_{p} \bigr\|_1 (x_p - x_{p-1})
\\
&=
 \int_{0}^{1}\! \|\pi(x) - \tilde{\pi}(x)\|_1 \, dx.
\end{align*}
This finishes the proof.
\qed

\medskip
\noindent
\textbf{Lipschitz continuity II.}
Next, let us consider the functional in (\ref{funp2}),
\begin{equation}
f_N^2(\pi) = \frac{1}{N}\smsp \e\log \sum_{\alpha\in\Natural^r} v_{\alpha} 
\exp \beta \sqrt{N} Y^\alpha.
\label{funp22}
\end{equation}
It actually does not depend on $N$ and, by (\ref{rearrange}), it can be represented as
\begin{equation}
f_N^2(\pi)
=
- \frac{\beta^2}{2}\sum_{k\leq \kappa} d_{k}^2 + \frac{\beta^2}{2}\int_{0}^1\! \|\pi(x)\|_{HS}^2\, dx.
\label{rearr2}
\end{equation}
In particular, it obviously satisfies
\begin{equation}
\bigl|f_N^2(\pi)-f_N^2(\tilde{\pi}) \bigr|
\leq
L \int_{0}^{1}\! \|\pi(x) - \tilde{\pi}(x)\|_1 \, dx.
\label{Lipschitz2}
\end{equation}
The equations (\ref{Lipschitz1}) and (\ref{Lipschitz2}) prove that these two types of functionals are uniformly Lipschitz on the set of discrete paths in $\Pi_d$ with respect to the metric $\Delta$ in (\ref{metric1}).

\medskip
\noindent
\textbf{Discretization with respect to $\Delta$.}
To conclude that these functionals can be extended to Lipschitz functionals on the entire $\Pi_d$, we need to observe that any path $\pi\in \Pi_d$ can be approximated by a discrete path with respect to $\Delta$. To see this, notice that for any $\gamma\in\Gamma_\kappa,$ $\|\gamma\|_1 \leq \kappa \tr(\gamma)$, because $|\gamma_{k,k'}|\leq (\gamma_{k,k}+\gamma_{k',k'})/2.$ For $\pi\in \Pi_d,$ $\pi(x')-\pi(x)\in\Gamma_\kappa$ for $x\leq x'$ and, therefore, 
$$
\|\pi(x')-\pi(x)\|_1 \leq \kappa \tr\bigl(\pi(x')-\pi(x)\bigr) = \kappa\bigl(\tr(\pi(x'))-\tr(\pi(x))\bigr).
$$
This implies that for any $x,x'\in [0,1]$, $\|\pi(x')-\pi(x)\|_1 \leq \kappa | \tr(\pi(x'))-\tr(\pi(x))|.$ Therefore, if we consider any discretization of the path $\pi$,
\begin{equation}
\tilde{\pi}(x) := \pi(x_p^*) \,\mbox{ for }\, x_{p-1}< x\leq x_p, 0\leq p\leq r,
\label{pitilde}
\end{equation}
for arbitrary choice of points $(x_p^*)$ inside these intervals, then
\begin{equation}
\Delta(\pi, \tilde{\pi}) = \int_{0}^{1}\! \|\pi(x) - \tilde{\pi}(x)\|_1 \, dx
\leq \kappa  \int_{0}^{1}\!  \bigl| \tr(\pi(x))-\tr(\tilde{\pi}(x)) \bigr| \, dx.
\label{Deltapitilde}
\end{equation}
Since the function $\tr(\pi(x))$ is non-decreasing for $\pi\in \Pi_d$, we can, obviously, make the right hand side as small as we like by an appropriate choice of sequences $(x_p)$ and $(x_p^*).$

\medskip
\noindent
\textbf{Another type of continuity.} The functionals considered above will appear as the limit of some functionals defined on finite size systems in terms of some Gaussian processes whose covariance structure becomes related to the Ruelle probability cascades only in the limit, due to the key result that will be proved below. As in the classical Sherrington-Kirkpatrick model, the covariance will be a function of the overlaps and the functionals will be continuous with respect to the distribution of the overlaps, allowing us to express them in the limit in terms of the functionals considered above. 

Such continuity properties are quite standard, and here we will only remind their general form.  For example, for a fixed $N\geq 1$, consider a functional resembling (\ref{fNS}),
\begin{equation}
f_1 = \frac{1}{N} \e \log \sum_{\alpha\in \A} w_\alpha \sum_{\sigma\in S} \exp \beta \sum_{i\leq N} Z_{i}^\alpha(\sigma_i)
\label{fNabs}
\end{equation}
for an arbitrary non-empty subset $S\subseteq \{1,\ldots,\kappa\}^N$. Here, the random weights $(w_\alpha)_{\alpha\in\A}$ define some random probability distribution $G$ on a countable (or finite) set $\A$, and $G$ is independent of the Gaussian processes $Z_{i}^\alpha = (Z_{i}^\alpha(k))_{k\leq \kappa}$. These Gaussian processes are independent copies of some Gaussian process $Z^\alpha = (Z^\alpha(k))_{k\leq \kappa}$ with the covariance
\begin{equation}
\Cov(Z^{\alpha^1}, Z^{\alpha^2}) = C_z\bigl(R_{\alpha^1,\alpha^2}\bigr)
\label{CzCov}
\end{equation}
for some continuous functions $C_z$ of the `overlaps' $R_{\alpha^1,\alpha^2} = \bigl(R_{\alpha^1,\alpha^2}^{k,k'}\bigr)_{k,k'\leq \kappa}.$ The array 
\begin{equation}
R_\A=\bigl(R_{\alpha^1,\alpha^2}\bigr)_{\alpha^1,\alpha^2 \in \A}
\label{Ralpha}
\end{equation}
is non-random with bounded entries, and at this moment we think of it as corresponding to some abstract overlap structure, for example, some infinite Gram array. Similarly, we can define the analogue of the functional in (\ref{funp22}),
 \begin{equation}
f_2 = \frac{1}{N}\smsp \e\log \sum_{\alpha\in\A} w_{\alpha}  \exp \beta \sqrt{N} Y^\alpha
\label{fN2abs}
\end{equation}
corresponding to the covariance
\begin{equation}
\Cov(Y^{\alpha^1}, Y^{\alpha^2}) = C_y\bigl(R_{\alpha^1,\alpha^2}\bigr)
\label{CyCov}
\end{equation}
for some continuous functions $C_y$. Let $(\alpha(\ell))_{\ell\geq 1}$ be i.i.d. indices sampled from the distribution $G(\alpha)=w_\alpha$ on $\A$, and let
\begin{equation}
R=\bigl(R_{\alpha(\ell),\alpha(\ell')}\bigr)_{\ell,\ell'\geq 1}.
\label{Rarandom}
\end{equation}
The following holds.
\begin{lemma}\label{LemCont2}
The quantities $f_1$ in (\ref{fNabs}) and $f_2$ in (\ref{fN2abs}) are continuous functionals of the distribution of the array $R$ in (\ref{Rarandom}) under $\e G^{\otimes\infty}.$ In other words, they depend on $\A,$ $(w_\alpha)_{\alpha\in \A}$ and the covariance structure $R_\A$ in (\ref{Ralpha}) only through the distribution of the array $R$ in (\ref{Rarandom}) under $\e G^{\otimes\infty}.$
\end{lemma}
Of course, these functionals depend on $N$, the set $S$, $\beta$, $C_z$ in (\ref{CzCov}). More specifically, the lemma says that, for any $\eps>0$, there exists $n\geq 1$ and a continuous bounded function $f_{1,\eps}$ of the array $R^n = (R_{\alpha(\ell),\alpha(\ell')})_{1\leq \ell,\ell'\leq n}$ such that 
$$
|f_1-\e f_{1,\eps}(R^n)|\leq \eps.
$$
The function $f_{1,\eps}$ depends only on $\eps$, $N$, $S$, $\beta$, and $C_z$. The same holds for $f_2$. The proof is omitted as it is almost identical, for example, to the proof of Theorem 1.3 in \cite{SKmodel}.

\section{A new family of the Ghirlanda-Guerra identities} \label{Sec5label}

The proof of the lower bound in Theorem \ref{ThFE} will rely on the synchronization mechanism based on a new family of  identities of the Ghirlanda-Guerra type \cite{GG}. These resemble the multi-species identities in \cite{PMS}, but the difference is that now we deal with blocks of overlaps, and to study their matrix properties we need new type of identities. These identities arise via a small perturbation of the Hamiltonian that we will now define. 

For $p\geq 1$, we will use the following notation,
$$
e = (i_1,\ldots,i_p)\in \{1,\ldots, N\}^p,\,\,
\sigma_e = (\sigma_{i_1},\ldots,\sigma_{i_p})
$$
for a given $\sigma\in\{1,\ldots,\kappa\}^N.$ Given $\lambda\in\Reals^\kappa$, we will denote
$$
S_\lambda(\sigma_e) = \sum_{k\leq \kappa} \lambda_k \I(\sigma_{i_1}=k)\cdots\I(\sigma_{i_p}=k).
$$
Given $n\geq 0$ and $I=(e_1,\ldots, e_{n})\in (\{1,\ldots, N\}^p)^{n}$, we will denote
$$
S_\lambda(\sigma_I) = S_\lambda(\sigma_{e_1})\cdots S_\lambda(\sigma_{e_n}).
$$
For integer $m\geq 1$ and $n_1,\ldots,n_m\geq 1$, let $I_j = (e_1,\ldots, e_{n_j})\in (\{1,\ldots, N\}^p)^{n_j}$ and $\lambda^j\in\Reals^\kappa$ for $1\leq j\leq m$ and consider the following Hamiltonian
\begin{equation}
h_{\theta}(\sigma) = \frac{1}{N^{p(n_1+\ldots+n_m)/2}}
\sum_{I_1,\ldots,I_m} g_{I_1,\ldots,I_m}S_{\lambda^1}(\sigma_{I_1})\cdots S_{\lambda^m}(\sigma_{I_m}),
\label{htheta}
\end{equation}
where $g_{I_1,\ldots,I_m}$ are standard Gaussian random variables independent for different choices of the indices. For simplicity of notation, we denoted the list of all parameters of the Hamiltonian by
\begin{equation}
\theta = (p,m,n_1,\ldots,n_m,\lambda^1,\ldots,\lambda^m).
\label{thetapar}
\end{equation}
The covariance $C^\theta_{\ell,\ell'}: = \Cov(h_{\theta}(\sigma^\ell),h_{\theta}(\sigma^{\ell'})) $ of the Gaussian process (\ref{htheta}) equals
$$
C^\theta_{\ell,\ell'} = 
\prod_{j\leq m}\frac{1}{N^{p n_j}} 
\sum_{I_j} S_{\lambda^j}(\sigma_{I_j}^\ell) S_{\lambda^j}(\sigma_{I_j}^{\ell'})
=
\prod_{j\leq m}\Bigl(\frac{1}{N^p}\sum_{e} S_{\lambda^j}(\sigma_e^\ell) S_{\lambda^j}(\sigma_e^{\ell'}) \Bigr)^{n_j}.
$$
If we recall the notation for the matrix $R_{\ell,\ell'}$ in (\ref{Rll}) of overlaps (\ref{Rllkk}) then, for $\lambda\in\Reals^k$, we can rewrite
\begin{align*}
\frac{1}{N^p}\sum_{e} S_{\lambda}(\sigma_e^\ell) S_{\lambda}(\sigma_e^{\ell'}) 
&=
\sum_{k,k'\leq \kappa} \lambda_k\lambda_{k'}
\frac{1}{N^p}\sum_{i_1,\ldots,i_p} \prod_{r\leq p}
\I(\sigma_{i_r}^\ell=k)\I(\sigma_{i_r}^{\ell'}=k')
\\
&=
\sum_{k,k'\leq \kappa} \lambda_k\lambda_{k'} (R^{k,k'}_{\ell,\ell'})^p
=
\bigl( R_{\ell,\ell'}^{\circ p} \lambda,\lambda \bigr),
\end{align*}
where $A^{\circ p}$ denotes the Hadamard (element-wise) $p^{\mathrm{th}}$ power of the matrix $A$. Hence, the covariance can be written as
\begin{equation}
C^\theta_{\ell,\ell'}
=
\prod_{j\leq m}\bigl( R_{\ell,\ell'}^{\circ p} \lambda^j,\lambda^j \bigr)^{n_j}
\label{covhtheta}
\end{equation}
for any configurations $\sigma^\ell,\sigma^{\ell'}\in \{1,\ldots, \kappa\}^N.$

\medskip
\noindent
\emph{Definition.}
Let $\Theta$ be a collection of all $\theta$ of the type (\ref{thetapar}) with $p\geq 1$, $m\geq 1$, $n_1,\ldots,n_m\geq 1$, and $\lambda^1,\ldots,\lambda^m$ taking values in $([-1,1]\cap \mathbb{Q})^\kappa$ with all rational coordinates. 
\qed

\medskip
Let us consider a one-to-one function $j_0:([-1,1]\cap \mathbb{Q})^\kappa \to\Natural$ and let
$$
j(\theta) := p+n_1+\ldots+n_m+j_0(\lambda_1)+\ldots+j_0(\lambda_m)+4m.
$$
Let $(u_{\theta})_{\theta \in \Theta}$ be i.i.d. random variables uniform on the interval $[1,2]$ and define a Hamiltonian
\begin{equation}
h_{N}(\sigma) = \sum_{\theta\in\Theta} 2^{-j(\theta)} u_{\theta}\, h_{\theta}(\sigma).
\label{hNw}
\end{equation}
Conditionally on $u=(u_{\theta})_{\theta\in \Theta}$, this is a Gaussian process with the variance bounded by $1$. The Hamiltonian $h_N(\sigma)$ will be used as a perturbation of the model, which means that, instead of $H_N(\sigma)$ in (\ref{SKH}), from now on we will consider the perturbed Hamiltonian 
\begin{equation}
H_N^{\mathrm{pert}}(\sigma) = H_N(\sigma) + s_N h_N(\sigma),
\label{Hpert}
\end{equation} 
where $s_N=N^{\gamma}$ for any $1/4<\gamma<1/2$. With this choice, the perturbation term is negligible from the point of view of computation of the free energy because $\lim_{N\to\infty} N^{-1} s_N^2 = 0$.

As in the classical Sherrington-Kirkpatrick model, the perturbation term $h_N(\sigma)$ is introduced to ensure the validity of the Ghirlanda-Guerra identities \cite{GG} for the Gibbs measure. The main difference is that now we will work with the Gibbs measure restricted to the configurations with fixed state sizes. Recall the definition of the set $\Sigma(d)$ in (\ref{SigmaD}) and $\DD_N$ in (\ref{DDN}). Let us consider arbitrary $d^N \in \DD_N$ and define the Gibbs measure on $\Sigma(d^N)$ by
\begin{equation}
G_{d^N}(\sigma) = \frac{\exp\, H_N^{\mathrm{pert}} (\sigma)}{Z_N(d^N)}
\,\mbox{ where }\,
Z_N(d^N) = \sum_{\sigma\in\Sigma(d^N)} \exp\, H_N^{\mathrm{pert}}(\sigma),
\label{GNpert}
\end{equation}
of course, in the form adapted to the present model. As usual, we will denote the average with respect to $G_{d^N}^{\otimes \infty}$ by $\la\, \cdot\, \ra$. Now, given $n\geq 2,$ let (recall the definition of the matrices $R_{\ell,\ell'}$ in (\ref{Rll}))
$$
R^n = \bigl(R_{\ell, \ell'}\bigr)_{\ell,\ell'\leq n}
$$
and consider an arbitrary bounded measurable function $f=f(R^n)$. For $\theta \in \Theta$, let 
\begin{equation}
\varDelta(f,n,\theta) = 
\Bigl|
{\e}  \bigl\la f C_{1,n+1}^{\theta} \bigr\ra -  \frac{1}{n} {\e} \bigl\la f \bigr\ra {\e} \bigl\la C_{1,2}^{\theta} \bigr\ra - \frac{1}{n}\sum_{\ell=2}^{n}{\e} \bigl\la f C_{1,\ell}^{\theta} \bigr\ra
\Bigr|,
\label{GGfinite}
\end{equation}
where ${\e}$ denotes the expectation conditionally on the i.i.d. uniform sequence $u=(u_{\theta})_{\theta\in \Theta}$. If we denote by $\e_u$ the expectation with respect to $u$ then the following holds. 
\begin{lemma}\label{ThGG} 
For any $n\geq 2$ and any bounded measurable function $f=f(R^n)$, for all $\theta\in\Theta$,
\begin{equation}
\lim_{N\to\infty} \e_u \smsp \varDelta(f,n,\theta) = 0.
\label{GGxlim}
\end{equation} 
\end{lemma}
\textbf{Proof.} The proof is, essentially, identical to proof of Theorem 3.2 in \cite{SKmodel}. We only need to mention why restricting the Gibbs measure to the set of configurations $\Sigma(d^N)$ with fixed state sizes is so important. This is because the proof depends in a crucial way on the fact that the diagonal elements $C_{\ell,\ell}^{\theta}$ are constant independent of $\sigma^\ell.$ In our case,
$$
C_{\ell,\ell}^{\theta} = 
\prod_{j\leq m}\bigl( R_{\ell,\ell}^{\circ p} \lambda^j,\lambda^j \bigr)^{n_j}
=
\prod_{j\leq m}\Bigl( \diag\bigl((d_1^N)^p,\ldots, (d_\kappa^N)^p\bigr) \lambda^j,\lambda^j \Bigr)^{n_j}
$$
are, indeed, independent of the configuration $\sigma^\ell$ due to the constraint $\sigma\in\Sigma(d^N)$. Besides this observation, the rest of the argument is exactly the same and, for a given $\theta\in\Theta$, the equation (\ref{GGxlim}) is obtained by utilizing the term $h_{\theta}(\sigma)$ in the perturbation (\ref{hNw}). 
\qed

\medskip
Using (\ref{GGxlim}), one can choose a non-random sequence $u^N=(u^N_{\theta})_{\theta\in \Theta}\in [1,2]^\Theta$ such that 
\begin{equation}
\lim_{N\to\infty} \smsp \varDelta(f,n,\theta) = 0
\,\mbox{ for all }\, \theta\in\Theta
\label{GGxlim2}
\end{equation} 
for the Gibbs measure $G_N$ with the parameters $u$ in the perturbation (\ref{hNw}) equal to $u^N$ rather than random. In fact, the choice of $u^N$ will be made below in a special way to coordinate with the cavity computations in the lower bound. Right now we will consider any such sequence $u^N$. 

Let us now consider any subsequence $(N_k)_{k\geq 1}$ along which the array $(R_{\ell,\ell'})_{\ell,\ell'\geq 1}$ of the $\kappa\times\kappa$ overlap matrices in (\ref{Rll}) converges in distribution under the measure $\e G_N^{\otimes\infty}$. We will continue to use the same notation as in (\ref{Rll}) and (\ref{covhtheta}),
\begin{equation}
R_{\ell,\ell'} = (R_{\ell,\ell'}^{k,k'})_{k,k'\leq \kappa},\,
R^n = \bigl(R_{\ell, \ell'}\bigr)_{\ell,\ell'\leq n},\,
C^\theta_{\ell,\ell'} =
\prod_{j\leq m} \bigl( R_{\ell,\ell'}^{\circ p} \lambda^j,\lambda^j \bigr)^{n_j},
\label{Rwlim}
\end{equation}
for the limiting random array. Then the equations (\ref{GGfinite}) and (\ref{GGxlim2}) imply that \begin{equation}
\e f(R^n) C^\theta_{1,n+1}
=
 \frac{1}{n}\e f(R^n)  \e C^\theta_{1,2}
+ \frac{1}{n}\sum_{\ell=2}^{n}\e f(R^n) C^\theta_{1,\ell}
\label{GGwp}
\end{equation}
for all $\theta\in\Theta$. Since $C^\theta_{\ell,\ell'}$ is a continuous function of $\lambda^j\in [-1,1]^{\kappa}$ for $j\leq m$, (\ref{GGwp}) holds a posteriori for all values of $\lambda^j$, not only with rational coordinates. 

For any $p\geq 1$, $\lambda^1,\ldots,\lambda^m \in [-1,1]^{\kappa}$ and a bounded measurable function $\varphi\colon \Reals^m\to\Reals$, let
\begin{equation}
Q_{\ell,\ell'} := \varphi \Bigl( \bigl( R_{\ell,\ell'}^{\circ p} \lambda^1,\lambda^1 \bigr),\ldots, \bigl( R_{\ell,\ell'}^{\circ p} \lambda^m,\lambda^m \bigr) \Bigr).
\label{Qphi}
\end{equation} 
Then the following version of the Ghirlanda-Guerra identities holds.
\begin{theorem}\label{ThGGms}
For any $n\geq 2$ and any bounded measurable function $f=f(R^n)$,
\begin{equation}
\e f(R^n) Q_{1,n+1}
=
 \frac{1}{n}\e f(R^n)  \smsp \e Q_{1,2}
+ \frac{1}{n}\sum_{\ell=2}^{n}\e f(R^n) Q_{1,\ell}.
\label{GGms}
\end{equation}
\end{theorem}
\textbf{Proof.} By (\ref{GGwp}), this holds for all functions of the type
$
\varphi(y_1,\ldots,y_m)=y_1^{n_1}\cdots y_{m}^{n_m}.
$
Approximating continuous functions by polynomials, this implies (\ref{GGms}) for continuous functions $\varphi$ and the general case follows.
\qed

\section{Synchronizing the block of overlaps}\label{Sec6label}

In this section, we will prove the following result, which will be the main tool in the computation of the free energy. Recall the notation $\Gamma_\kappa$ in (\ref{GammaK}).
\begin{theorem}\label{Th2}
Suppose that the array $(R_{\ell,\ell'})_{\ell,\ell'\geq 1}$ in (\ref{Rwlim}) satisfies (\ref{GGms}) for all choices of parameters there. Then there exists a function $\Phi\colon \Reals^+ \to \Gamma_\kappa$ such that
\begin{equation}
R_{\ell,\ell'} = \Phi\bigl(\tr(R_{\ell,\ell'}) \bigr)
\end{equation}
almost surely. Moreover, one can take $\Phi$ to be non-decreasing in $\Gamma_\kappa$, 
$$
\Phi(x')-\Phi(x)\in \Gamma_\kappa
\,\mbox{ for all $x\leq x'$},
$$
and Lipschitz continuous,
$$
\|\Phi(x') - \Phi(x)\|_1 \leq L_\kappa |x'-x|
$$ 
for some constant $L_\kappa$ that depends only on $\kappa.$
\end{theorem}
The function $\Phi$ here is, of course, not universal and depends on the distribution of the array $(R_{\ell,\ell'}).$ 

\medskip
\noindent
\emph{Remark.} Notice that $0\leq R_{\ell,\ell'}^{k,k}\leq d_k$ because $R_{\ell,\ell}^{k,k}=R_{\ell',\ell'}^{k,k}=d_k$ and the entire overlap array is positive and positive-semidefinite. Therefore, $0\leq \tr(R_{\ell,\ell'})\leq \sum_{k\leq \kappa} d_k \leq 1$ and the function $\Phi$ can be defined on the interval $[0,1]$ instead of $\Reals^+.$

\medskip
The proof of Theorem \ref{Th2} will utilize the following results obtained in \cite{PMS}. For a fixed $p\geq 1$ and fixed $\lambda^1,\ldots,\lambda^m \in [-1,1]^\kappa$, let us consider the arrays
\begin{equation}
\Lambda_{\ell,\ell'}^j = \bigl( R_{\ell,\ell'}^{\circ p} \lambda^j,\lambda^j \bigr),\,
\Lambda_{\ell,\ell'} = \sum_{j\leq m} \Lambda_{\ell,\ell'}^j
\label{Aoverlap}
\end{equation}
indexed by $\ell,\ell'\geq 1.$ All of these arrays are symmetric, positive-semidefinite, and exchangeable in the sense that their distribution is invariant under the same permutation of finitely many rows and columns. All of these properties hold trivially before we pass to the limit $N_k\to\infty$ (see paragraph above (\ref{Rwlim})) and are inherited by the limiting array. The proof of Theorem \ref{Th2} will use some results for such arrays from \cite{PMS} (Theorem 4 and Lemma 2 there), which will be summarized in the next lemma. We omit the proof, because it can be carried over to the present setting without any modifications.
\begin{lemma}\label{ThSynch}
Suppose that the array $(R_{\ell,\ell'})$ satisfies (\ref{GGms}). Then the following hold.
\begin{enumerate}
\item[(i)] 
With probability one, if $\Lambda^j_{\ell,\ell'}> \Lambda^j_{\ell,\ell''}$ for some $j\leq m$ then $\Lambda^j_{\ell,\ell'}\geq \Lambda^j_{\ell,\ell''}$ for all $j\leq m$.

\item[(ii)] There exist non-decreasing $1$-Lipschitz functions $L_j\colon \Reals^+\to \Reals^+$ for $j\leq m$ such that, with probability one, $\Lambda^j_{\ell,\ell'} = L_j(\Lambda_{\ell,\ell'})$. 
\end{enumerate}
\end{lemma}
The reason we can consider the domain and range of $L_j$ to be $\Reals^+$ is because, by (\ref{GGms}), each array in (\ref{Aoverlap}) by itself satisfies the canonical Ghirlanda-Guerra identities \cite{GG} and, therefore, its entries are nonnegative by Talagrand's positivity principle (see Theorem 2.16 in \cite{SKmodel}). Equipped with Lemma \ref{ThSynch}, we proceed to prove Theorem \ref{Th2}.

\medskip
\noindent
\textbf{Proof of Theorem \ref{Th2}.} \emph{Step 1.} In (\ref{Aoverlap}), let us take $p=1$, $m=\kappa$ and let $\lambda^1=e_1,\ldots, \lambda^\kappa=e_\kappa$ be the standard basis in $\Reals^\kappa.$ With this choice of parameters, 
$$
\Lambda^k_{\ell,\ell'} = R_{\ell,\ell'}^{k,k}
\,\mbox{ and }\,
\Lambda_{\ell,\ell'} = \sum_{k\leq \kappa} R_{\ell,\ell'}^{k,k} = \tr(R_{\ell,\ell'}).
$$
Lemma \ref{ThSynch} (ii) implies that there exist non-decreasing $1$-Lipschitz functions $L_k$ such that 
\begin{equation}
R_{\ell,\ell'}^{k,k} = L_k\bigl(\tr(R_{\ell,\ell'})\bigr)
\label{Step1}
\end{equation}
with probability one.

\emph{Step 2.} Let us fix two indices $k,k'\leq \kappa.$ In (\ref{Aoverlap}), let us take $p=1$, $m=2$ and let $\lambda^1 = e_k + e_{k'}$ and $\lambda^2=e_k-e_{k'}.$ Then
\begin{align*}
\Lambda^1_{\ell,\ell'} &= R_{\ell,\ell'}^{k,k}+R_{\ell,\ell'}^{k',k'}+R_{\ell,\ell'}^{k,k'}+R_{\ell,\ell'}^{k',k},
\\
\Lambda^2_{\ell,\ell'} &= R_{\ell,\ell'}^{k,k}+R_{\ell,\ell'}^{k',k'}-R_{\ell,\ell'}^{k,k'}-R_{\ell,\ell'}^{k',k},
\\
\Lambda_{\ell,\ell'} &= 2(R_{\ell,\ell'}^{k,k}+R_{\ell,\ell'}^{k',k'}) = L\bigl(\tr(R_{\ell,\ell'})\bigr),
\end{align*}
where the function $L=2(L_k+L_{k'})$ and the last equality follows from (\ref{Step1}). Since 
$$
R_{\ell,\ell'}^{k,k'}+R_{\ell,\ell'}^{k',k} = \Lambda^1_{\ell,\ell'} - \frac{1}{2}\Lambda_{\ell,\ell'},
$$
Lemma \ref{ThSynch} (ii) implies that there exist a Lipschitz functions $L_{k,k'}$ (maybe, not monotone) such that 
\begin{equation}
R_{\ell,\ell'}^{k,k'}+R_{\ell,\ell'}^{k',k} = L_{k,k'}\bigl(\tr(R_{\ell,\ell'})\bigr)
\label{Step2}
\end{equation}
with probability one.

\emph{Step 3.} If in the above two steps we take $p=2$ then the same arguments shows that, for any $k,k'\leq \kappa$, 
$
(R_{\ell,\ell'}^{k,k'})^2+(R_{\ell,\ell'}^{k',k})^2
$
is a Lipschitz function of $\tr(R_{\ell,\ell'}^{\circ 2}).$ However, since each $R_{\ell,\ell'}^{k,k}$ is bounded and itself is a Lipschitz function of $\tr(R_{\ell,\ell'}),$ the trace $\tr(R_{\ell,\ell'}^{\circ 2})$ is also a Lipschitz function of $\tr(R_{\ell,\ell'}).$ Therefore, there exist a Lipschitz functions $L_{k,k'}'$ such that, with probability one,
\begin{equation}
(R_{\ell,\ell'}^{k,k'})^2+(R_{\ell,\ell'}^{k',k})^2 = L_{k,k'}'\bigl(\tr(R_{\ell,\ell'})\bigr).
\label{Step3}
\end{equation}

\emph{Step 4. } The systems of equations (\ref{Step2}) and (\ref{Step3}) is of the form $x+y=a, x^2+y^2=b$ and can be solved to find $\{x,y\}$ in terms of $a$ and $b$, 
$$
x,y = \frac{a\pm \sqrt{2b-a^2}}{2}.
$$
In other words, there exist two continuous functions $f_1, f_2$ (maybe, not Lipschitz) such that
\begin{equation}
\bigl\{R_{\ell,\ell'}^{k,k'}, R_{\ell,\ell'}^{k',k} \bigr\} = 
\bigl\{ f_1\bigl(\tr(R_{\ell,\ell'})\bigr), f_2\bigl(\tr(R_{\ell,\ell'})\bigr) \bigr\}.
\label{f1f2}
\end{equation}
(The functions, of course, depend on the indices $k,k'$.) The main obstacle in the proof is that we do not know which of these two overlaps takes which of the two values, so this does not quite reconstruct the matrix in terms of its trace $\tr(R_{\ell,\ell'})$. However, in the next step we will show that one can take $f_1=f_2$ and, in particular, by (\ref{Step2}),
\begin{equation}
R_{\ell,\ell'}^{k,k'} = R_{\ell,\ell'}^{k',k} = \frac{1}{2}L_{k,k'}\bigl(\tr(R_{\ell,\ell'})\bigr)
\label{Step4}
\end{equation}
with probability one.

\emph{Step 5. } The array $(\tr(R_{\ell,\ell'}))_{\ell,\ell'\geq 1}$ is symmetric, positive-semidefinite, exchangeable and, by (\ref{GGms}), satisfies the canonical Ghirlanda-Guerra identities \cite{GG}. Therefore, the results in Sections 2.4 in \cite{SKmodel} imply that it can be generated by the Ruelle probability cascades and, in particular, it satisfies what was called the duplication property in Sections 2.5 in \cite{SKmodel}. This means that if the support (of the distribution) of $\tr(R_{1,2})$ contains a point $q$ then support of the array $(\tr(R_{\ell,\ell'}))_{1\leq \ell<\ell' \leq n}$ contains the array with all entries equal to $q$ for any $n\geq 2.$ Recalling (\ref{Step1}) and (\ref{f1f2}), let us denote
$$
a = L_k(q), d = L_{k'}(q), b = f_1(q), c=f_2(q).
$$
By the above steps, the $2n\times 2n$ array consisting of $2\times 2$ blocks indexed by $1\leq \ell,\ell' \leq n$,
$$
\begin{bmatrix}
R_{\ell,\ell'}^{k,k} & R_{\ell,\ell'}^{k,k'} \\
R_{\ell,\ell'}^{k',k}  & R_{\ell,\ell'}^{k',k'} 
\end{bmatrix},
$$
will have in its support a $2n\times 2n$ array $A$ indexed by $1\leq \ell,\ell' \leq n$ consisting of $2\times 2$ blocks 
$$
a_{\ell,\ell'}=
\begin{bmatrix}
a_{\ell,\ell'}^{1,1} & a_{\ell,\ell'}^{1,2} \\
a_{\ell,\ell'}^{2,1} & a_{\ell,\ell'}^{2,2}
\end{bmatrix},
$$
where each of the blocks $a_{\ell,\ell'}$ is either
$$
\begin{bmatrix}
d_k & 0 \\
0 & d_{k'}
\end{bmatrix},\,
\begin{bmatrix}
a & b \\
c & d
\end{bmatrix}
\,\mbox{ or }\,
\begin{bmatrix}
a & c \\
b & d
\end{bmatrix},
$$ 
where the first choice corresponds to the diagonal blocks for $\ell=\ell'$, and the second and third correspond to $\ell\not = \ell'.$ Let us define a directed complete graph $(V,E)$ on $n$ vertices such that, for each pair $\ell\not =\ell'$, the edge is oriented $\ell\to\ell'$ or $\ell'\to\ell$ depending on whether
$$
a_{\ell,\ell'}=
\begin{bmatrix}
a & b \\
c & d
\end{bmatrix}
\,\mbox{ or }\,
a_{\ell,\ell'}=
\begin{bmatrix}
a & c \\
b & d
\end{bmatrix}.
$$

Directed complete graphs are called tournaments, and we will need to use the fact that, for large $n$, one can find two large disjoint subsets of vertices $V_1, V_2\subseteq V$ such that all edges between them are oriented from $V_1$ to $V_2.$  For example, we can use Theorem 3 in \cite{Erdos}, which states the following. Given two tournaments $S$ and $T$, $T$ is called $S$-free if $S$ is not a sub-tournament of $T.$ Theorem 3 in \cite{Erdos} states that if $T$ is $S$-free, $\mathrm{card}(S)=m$ and $\mathrm{card}(T)=n$ then one can find two disjoint subsets $V_1, V_2$ in $T$ of cardinality 
$$
\mathrm{card}(V_1) = \mathrm{card}(V_2) =
\Bigl \lfloor{(n/m)^{1/(m-1)}} \Bigr\rfloor
$$
with all edges between them oriented from $V_1$ to $V_2.$  Take $S$ consisting of two groups of cardinality $m$ with edges between groups oriented in one direction, and in arbitrary fashion within groups. If $T$ contains $S$ then is has two desired subsets $V_1, V_2$ of cardinality $m$. If not then, by the above claim, it contains two such subsets of cardinality $\lfloor{(n/2m)^{1/(2m-1)}}\rfloor.$ If we set $(n/2m)^{1/(2m-1)}= m/2$, we find that $m$ is of order $\log n / \log\log n$, which means that, for large $n$, we can always find two large disjoint subsets with edges between them oriented in the same direction.

Let $V_1,V_2\subseteq V=\{1,\ldots, n\}$ be two such disjoint subsets of size $m$ in the above graph. This means that for all $\ell\in V_1$ and $\ell'\in V_2$,
$$
a_{\ell,\ell'}=
\begin{bmatrix}
a & b \\
c & d
\end{bmatrix}.
$$
On the other hand, recall that the entire $2n\times 2n$ array $A$ is positive-semidefinite since it belongs to the support of a positive-semidefinite random array. Therefore, we can find pairs of vectors $u_\ell,w_{\ell}$ for $1\leq \ell\leq n$ in a Hilbert space $H$ such that, for all $\ell,\ell'\leq n$,
$$
\begin{bmatrix}
(u_{\ell},u_{\ell'}) & (u_{\ell},w_{\ell'}) \\
(w_{\ell},u_{\ell'}) & (w_{\ell},w_{\ell'})
\end{bmatrix}
=
\begin{bmatrix}
a_{\ell,\ell'}^{1,1} & a_{\ell,\ell'}^{1,2} \\
a_{\ell,\ell'}^{2,1} & a_{\ell,\ell'}^{2,2}
\end{bmatrix}.
$$
In particular, for all $\ell,\ell'\leq n$,
\begin{equation}
(u_{\ell},u_{\ell}) = d_k, (w_{\ell},w_{\ell}) = d_{k'}, (u_{\ell},u_{\ell'}) = a, (w_{\ell},w_{\ell'}) = d,
\label{Step5eq1}
\end{equation}
and, by construction of the sets $V_1, V_2$, for all $\ell\in V_1$ and $\ell'\in V_2$,
\begin{equation}
\begin{bmatrix}
(u_{\ell},u_{\ell'}) & (u_{\ell},w_{\ell'}) \\
(w_{\ell},u_{\ell'}) & (w_{\ell},w_{\ell'})
\end{bmatrix}
=
\begin{bmatrix}
a & b \\
c & d
\end{bmatrix}.
\label{Step5eq2}
\end{equation}
For $j=1,2$, let us consider the barycenters of these collections of vectors
$$
U_j = \frac{1}{m}\sum_{\ell\in V_j} u_\ell
\,\mbox{ and }\,
W_j = \frac{1}{m}\sum_{\ell\in V_j} w_\ell. 
$$
On the one hand, using the equation (\ref{Step5eq2}),
$$
\begin{bmatrix}
(U_1,U_2) & (U_1,W_2) \\
(W_1,U_2) & (W_1,W_2)
\end{bmatrix}
=
\begin{bmatrix}
a & b \\
c & d
\end{bmatrix}.
$$ 
On the other hand, using (\ref{Step5eq1}),
$$
\|U_1 - U_2\|^2 = \frac{1}{m^2}\Bigl\|\sum_{\ell\in V_1} u_\ell-\sum_{\ell\in V_2} u_\ell\Bigr\|^2
=
\frac{2(d_k+a)}{m}
$$
and
$$
\|W_1 - W_2\|^2 = \frac{1}{m^2}\Bigl\|\sum_{\ell\in V_1} w_\ell-\sum_{\ell\in V_2} w_\ell\Bigr\|^2
=
\frac{2(d_{k'}+d)}{m}.
$$
For large $m$, this implies that $U_1\approx U_2$ and $W_1\approx W_2$ and, therefore, $(U_1,W_2) \approx (W_1,U_2).$ Letting $m\to \infty$ proves that $b=c$, so we can take $f_1=f_2$ in (\ref{f1f2}), also proving (\ref{Step4}).

\emph{Step 6.} We proved that there exists a Lipschitz function $\Phi$ on $\Reals^+$ with values in the set of symmetric $\kappa\times\kappa$ matrices such that
\begin{equation}
R_{\ell,\ell'} = \Phi\bigl(\tr(R_{\ell,\ell'}) \bigr)
\end{equation}
almost surely. It remains to show that it is also non-decreasing in the space of Gram matrices (\ref{GammaK}), $\Phi(x')-\Phi(x)\in \Gamma_\kappa$ for all $x\leq x'$. Let us first prove this for $x, x'$ in the support of the distribution of $\tr(R_{1,2}).$ Suppose that $x< x'$ but $\Phi(x')-\Phi(x)\not \in \Gamma_\kappa$. Then, there exists $\lambda\in [-1,1]^\kappa$, such that 
$$
(\Phi(x)\lambda,\lambda)>(\Phi(x')\lambda,\lambda).
$$
In (\ref{Aoverlap}), let us take $p=1$, $m=\kappa+1$ and let $\lambda^1=e_1,\ldots, \lambda^\kappa=e_\kappa$ be the standard basis in $\Reals^\kappa$ and $\lambda^{\kappa+1} = \lambda$. With this choice of parameters, 
$$
\Lambda^k_{\ell,\ell'} = R_{\ell,\ell'}^{k,k} \,\mbox{ for }\, k\leq \kappa,
\Lambda^{\kappa+1}_{\ell,\ell'} = (R_{\ell,\ell'}\lambda,\lambda),
\,\mbox{ and }\,
\Lambda_{\ell,\ell'} =\tr(R_{\ell,\ell'}) + (R_{\ell,\ell'}\lambda,\lambda).
$$
Recall that Step 5 started with the statement that the array $(\tr(R_{\ell,\ell'}))_{\ell,\ell'\geq 1}$ is symmetric, positive-semidefinite, exchangeable and, by (\ref{GGms}), satisfies the canonical Ghirlanda-Guerra identities \cite{GG}. As a result, it satisfied the duplication property. Another consequence of the Ghirlanda-Guerra identities from Lemma 2.7 in \cite{SKmodel} (this was first observed in \cite{PT}) states that, with probability one, the set $\{\tr(R_{1,\ell}) \mid \ell\geq 2 \}$ is a dense subset of the support of the distribution of $\tr(R_{1,2})$. This means that, for any $\eps>0$, we can find $\ell,\ell'\geq 2$ such that 
$$
|\tr(R_{1,\ell}) - x|\leq \eps \,\mbox{ and }\,  |\tr(R_{1,\ell'}) - x'|\leq \eps.
$$ 
For small enough $\eps$, this implies that $\tr(R_{1,\ell})<\tr(R_{1,\ell'})$. By Lemma \ref{ThSynch} (i), $R_{1,\ell}^{k,k}<R_{1,\ell'}^{k,k}$ for at least one $k\leq \kappa$ and again, by Lemma \ref{ThSynch} (i), 
$$
\Lambda^{\kappa+1}_{1,\ell} = (R_{1,\ell}\lambda,\lambda)
= \bigl(\Phi(\tr(R_{1,\ell}))\lambda,\lambda \bigr)
\leq
\Lambda^{\kappa+1}_{1,\ell'} = (R_{1,\ell'}\lambda,\lambda)
= \bigl(\Phi(\tr(R_{1,\ell'}))\lambda,\lambda \bigr).
$$
Letting $\eps\downarrow 0$, we get that $$(\Phi(x)\lambda,\lambda)\leq (\Phi(x')\lambda,\lambda),$$ contradicting the above assumption. This proves that $\Phi$ is non-decreasing on the support of the distribution of $\tr(R_{1,2})$. On each interval $(x,x')$ outside of the support with $x,x'$ in the support, we extend $\Phi$ by a linear interpolation between the values $\Phi(x)$ and $\Phi(x')$, which does not affect the monotonicity and Lipschitz properties. This finishes the proof.
\qed

\section{Lower bound via cavity computations}\label{Sec7label}

Finally, to prove the matching lower bound we will combine the structural results for the overlaps proved above with standard cavity computations. We will start with an obvious inequality $F_N \geq F_N(d^N)$, for any $d^N\in\DD_N.$ For a fixed $d\in\DD$, we will choose $d^N$ converging to $d$ in such a way that $d^N_k =0$ whenever $d_k=0$ and, otherwise,
\begin{equation}
|d^N_k - d_k| \leq \frac{L_\kappa}{N}
\,\mbox{ for all }\, k\leq \kappa
\label{near}
\end{equation}
for some constant $L_\kappa$ that depends only on $\kappa.$ The next step is to use the inequality,
\begin{equation}
\liminf_{N\to\infty} F_N(d^N) \geq \frac{1}{M}\liminf_{N\to\infty} \Bigl(\e \log Z_{N+M}(d^{N+M}) - \e \log Z_{N}(d^N)\Bigr),
\label{FNAN}
\end{equation}
where $M$ on the right hand side is fixed and where, for any $N$ and $d\in \DD_N,$ we denoted
$$
Z_{N}(d) = \sum_{\sigma\in\Sigma_N(d)} \exp\beta H_N(\sigma),
$$
where now we will make the dependence of $\Sigma(d)=\Sigma_N(d)$ on the dimension explicit. The infimum on the right hand side of (\ref{FNAN}) is achieved along some subsequence $(N_k)_{k\geq 1}$, but to simplify the notation we will keep writing $N$. Next, from this subsequence, for any fixed $M$, we will choose another subsequence as follows (slightly modifying Lemma 5 in \cite{SKcoupled}).
\begin{lemma}
There exists a sequence $\delta^M \in \DD_M$ such that 
\begin{equation}
|\delta^M_k - d_k| \leq \frac{2L_\kappa}{M}
\,\mbox{ for all }\, k\leq \kappa
\label{near2}
\end{equation}
and such that, for each $M\geq 1,$
\begin{equation}
N d^{N} + M\delta^M = (N+M)d^{N+M}
\label{uprime}
\end{equation}
for infinitely many $N\geq 1.$ 
\end{lemma}
{\bf Proof.}
For a fixed $M,$ consider a sequence $\delta^M(N)$ defined by (\ref{uprime}),
$$
M\delta^M(N) := (N+M)d^{N+M} -N d^{N}.
$$
Subtracting $Md$ on both sides, for all $k\leq \kappa$,
$$
M(\delta^M(N)_k- d_k) = (N+M)(d^{N+M}_k - d_k) - N(d^{N}_k - d_k)
$$
and, therefore, (\ref{near}) implies that $M |\delta^M(N)_k- d_k)| \leq 2L_\kappa.$ For a fixed $M\geq 1$, this implies that the sequence $(M\delta^M(N))_{N\geq 1}$ takes a finite number of values and we can find infinitely many $N$ with the same value, denoted $M\delta^M$. Notice that, by construction, $\delta^M_k = 0$ whenever $d_k=0$ and, otherwise, $\delta^M_k \geq d_k - 2M^{-1}L_\kappa >0$ for $M$ large enough.
\qed

\medskip
For a fixed $M$, let us take a subsequence of $N$ found in (\ref{uprime}) and, again, for simplicity of notation we will keep writing $N$. The equation (\ref{uprime}) implies that
$$
\Sigma_{N+M}(d^{N+M}) \supseteq \Sigma_N(d^N)\times \Sigma_M(\delta^M).
$$
If we represent configurations $\rho\in \{1,\ldots,\kappa\}^{N+M}$ as $\rho=(\sigma,\eps)$ for $\sigma\in\{1,\ldots,\kappa\}^{N}$ and $\eps \in \{1,\ldots,\kappa\}^{M}$, this inclusion implies that
$$
Z_{N+M}(d^{N+M}) \geq \sum_{\sigma\in\Sigma_N(d^N)} \sum_{\eps\in\Sigma_M(\delta^M)} \exp\beta H_{N+M}(\sigma,\eps).
$$
This inequality is in the right direction for the purpose of decreasing the lower bound in (\ref{FNAN}) to
\begin{equation}
\lim_{N\to\infty} \frac{1}{M} \Bigl(
\e \log \sum_{\sigma\in\Sigma_N(d^N)} \sum_{\eps\in\Sigma_M(\delta^M)} \exp\beta H_{N+M}(\sigma,\eps)
- \e \log \sum_{\sigma\in\Sigma_N(d^N)} \exp\beta H_{N}(\sigma)
\Bigr),
\label{FMNAN}
\end{equation}
where, for a fixed $M$, the limit here is over some subsequence determined above.

Next, we separate a common part in the Hamiltonians $H_{N+M}(\sigma,\eps)$ and $H_{N}(\sigma)$ as in the usual Aizenman-Sim-Starr representation \cite{AS2} (see e.g. Section 1.3 in \cite{SKmodel}). First, we can separate
\begin{equation}
H_{N+M}(\sigma,\eps) = H_N'(\sigma) + \sum_{i\leq M} Z_{i}^{\sigma}(\eps_i) + r(\eps)
\label{decomp1}
\end{equation}
into three types of terms, 
\begin{align}
H_N'(\sigma) &= \frac{1}{\sqrt{N+M}} \sum_{1\leq i,j\leq N} g_{ij}\I(\sigma_i =\sigma_j),
\label{commonH}
\\
Z_{i}^{\sigma}(\eps_i) &=  \frac{1}{\sqrt{N+M}} \sum_{1\leq j\leq N} \bigl(g_{N+i,j} + g_{j,N+i} \bigr) \I(\sigma_j =\eps_i)
\nonumber
\\
r(\eps) &= \frac{1}{\sqrt{N+M}} \sum_{1\leq i,j\leq M} g_{N+i,N+j}\I(\eps_i = \eps_j).
\nonumber
\end{align}
One the other hand, the Gaussian process $H_N(\sigma)$ on $\{1,\ldots,\kappa\}^N$ can be decomposed into a sum of two independent Gaussian processes (in distribution),
\begin{equation}
H_N(\sigma) \stackrel{d}{=}
H_N'(\sigma) + \sqrt{M} Y^{\sigma},
\label{commonH2}
\end{equation}
where $H_N'(\sigma)$ was defined in (\ref{commonH}) and
\begin{equation}
Y^{\sigma} = 
 \frac{1}{\sqrt{N(N+M)}} \sum_{1\leq i,j\leq N} g_{ij}'\I(\sigma_i =\sigma_j),
 \end{equation}
where $(g_{ij}')$ are independent copies of the Gaussian random variables $(g_{ij})$. The term $r(\eps)$ can be omitted because it is of a small order as $N\to\infty$.

Consider the Gibbs measure on $\Sigma_N(d^N)$ corresponding to the Hamiltonian $H_N'(\sigma)$ in (\ref{commonH}),
\begin{equation}
G_N'(\sigma) = \frac{\exp \beta H_N'(\sigma)}{Z_N'(d^N)},
\,\mbox{ where }\,
Z_N'(d^N) = \sum_{\sigma\in\Sigma_{N}(d^N)} \exp \beta H_N'(\sigma),
\label{MeasureGNprime}
\end{equation}
and let us denote by $\la\,\cdot\,\ra_N'$ the average with respect to $G_N'$. Using representation (\ref{decomp1}) and (\ref{commonH2}) and dividing inside both $\log$arithms by $Z_N'(d^N)$, we can rewrite (\ref{FMNAN}) as
\begin{equation}
\lim_{N\to\infty}
\frac{1}{M}
\Bigl(
\e \log \Bigl\la
\sum_{\eps\in\Sigma_M(\delta^M)} \exp\beta \sum_{i\leq M} Z_{i}^{\sigma}(\eps_i)
\Bigr\ra_N'
-
\e \log \Bigl\la \exp \beta \sqrt{M} Y^{\sigma} \Bigr\ra_N'
\Bigr).
\label{AS2repr}
\end{equation}
Both terms here are exactly of the form considered in Lemma \ref{LemCont2} above with $\alpha=\sigma$,
$\A=\Sigma_N(d^N)$ and $w_\alpha = G_N'(\sigma)$ there. Therefore, they can be viewed as continuous functionals of the overlap arrays determined by the covariance structure of the Gaussian processes $Z_i^\sigma=(Z_i^\sigma(k))_{k\leq\kappa}$ and $Y^\sigma$, which we now compute. 

First of all, for $k,k'\leq \kappa$,
\begin{equation}
\e Z_i^{\sigma^\ell}(k) Z_i^{\sigma^{\ell'}}(k')
= 
\frac{2}{N+M}\sum_{j\leq N} \I(\sigma^\ell_j = k) \I(\sigma^{\ell'}_j = k')
= 
2R_{\ell,\ell'}^{k,k'} + \bigO{N^{-1}}
\label{Covz}
\end{equation}
and, similarly to the computation in (\ref{CovHNsum}),
\begin{equation}
\e Y^{\sigma^\ell} Y^{\sigma^{\ell'}} = 
\frac{N}{N+M}
\sum_{k,k'\leq \kappa}(R_{\ell,\ell'}^{k,k'})^2
=
\sum_{k,k'\leq \kappa}(R_{\ell,\ell'}^{k,k'})^2
+ \bigO{N^{-1}}.
\label{Covy}
\end{equation}
These resemble the definition in (\ref{CD}) and, of course, one can redefine the processes $Z_i^\sigma$ and $Y^\sigma$ to have covariances without the lower order terms $\bigO{N^{-1}}$, which we now assume.

The same computation can be carried out just as easily in the case when the constrained free energy $F_N(d^N)$ in (\ref{FNAN}) corresponds to the perturbed Hamiltonian $H_N^{\mathrm{pert}}(\sigma)$ in (\ref{Hpert}) instead of the original Hamiltonian $H_N(\sigma)$. Moreover, since the perturbation term $s_N h_N(\sigma)$ in (\ref{Hpert}) is of smaller order, one can show that the perturbation term $s_{N+M} h_{N+M}(\sigma,\eps)$ that would appear in the first term in (\ref{FMNAN}) can be replaced by $s_N h_N(\sigma)$ and this only introduces some small order correction. All of this is standard and is explained, for example, in Section 3.5 in \cite{SKmodel}. In other words, if the Gibbs measure $G_N'$ in (\ref{MeasureGNprime}) corresponds to the perturbed Hamiltonian
$$
 H_N'(\sigma) + s_N h_N(\sigma)
$$
then the representation in (\ref{AS2repr}) still gives a lower bound on the constrained free energy along some subsequence. Also, in this case the expectation $\e$ in (\ref{AS2repr}) includes the average $\e_u$ in the uniform random variables $u=(u_{\theta})_{\theta\in\Theta}$ in the definition of the perturbation Hamiltonian (\ref{hNw}).

The proof of Lemma \ref{ThGG} applies verbatim to the measure $G_N'$ and right below the proof of Lemma \ref{ThGG} we mentioned that one can choose a non-random sequence $u^N=(u^N_{\theta})_{\theta\in \Theta}$ changing with $N$ such that (\ref{GGxlim2}) holds for the Gibbs measure $G_N'$ with the parameters $u$ in the perturbation Hamiltonian (\ref{hNw}) equal to $u^N$ rather than random. By Lemma 3.3 in \cite{SKmodel}, one can choose this sequence $u^N$ in such a way that the limit in (\ref{AS2repr}) is also not affected by fixing $u=u^N$ instead of averaging in $u$. 

Passing to another subsequence, we can assume that the distribution of the array $(R_{\ell,\ell'})_{\ell,\ell'\geq 1}$ in (\ref{Rll}) under $\e G_N^{\prime \otimes\infty}$ converges, and we will denote the array with this limiting distribution by $(R^{\infty}_{\ell,\ell'})_{\ell,\ell'\geq 1}$. By construction, this array satisfies the generalized Ghirlanda-Guerra identities in Theorem \ref{ThGGms} and, as a consequence, its structure can be described as in Theorem \ref{Th2}. Namely, there exists a function $\Phi\colon [0,1] \to \Gamma_\kappa$ such that
\begin{equation}
R^{\infty}_{\ell,\ell'} = \Phi\bigl(\tr(R^{\infty}_{\ell,\ell'}) \bigr)
\label{RLS}
\end{equation}
almost surely. The function $\Phi$ is non-decreasing, $\Phi(x')-\Phi(x)\in \Gamma_\kappa$ for all $x\leq x'$,  and Lipschitz, $\|\Phi(x') - \Phi(x)\|_1 \leq L_\kappa |x'-x|$ for some constant $L_\kappa$ that depends only on $\kappa.$

Let us denote the distribution function of $\tr(R^{\infty}_{1,2})$ by
\begin{equation}
\mu_\infty(q) = \p\bigl(\tr(R^{\infty}_{1,2})\leq q \bigr)
\label{zetainfty}
\end{equation}
and let $\mu_\infty^{-1}\colon [0,1] \to \Reals^+$ be its quantile transformation. Define
\begin{equation}
\pi_\infty(x):= \Phi\bigl(\mu_\infty^{-1}(x)\bigr),
\label{piinfty}
\end{equation}
which is an element of the family of paths $\Pi_d$ defined in (\ref{Pid}). Notice that, by (\ref{RLS}), 
$$
\tr\bigl(\pi_\infty(x)\bigr) = \mu_\infty^{-1}(x).
$$
Let us consider two sequences,
\begin{align}
x_{-1}= 0 &< x_0 < \ldots < x_{r-1} < x_r = 1,
\label{xstrictagain}
\\
0 &=q_0<  \ldots < q_{r-1}< q_r = 1,
\nonumber
\end{align}
such that $q_p=\mu_\infty^{-1}(x_p)$. Consider the distribution function $\mu$ defined by
\begin{equation}
\mu(q) = x_{p}
\,\mbox{ for }\,
q_p\leq q<q_{p+1}
\label{zetafop}
\end{equation}
and, similarly to (\ref{piinfty}), we define the corresponding path in $\Pi_d$ by
\begin{equation}
\pi(x):= \Phi\bigl(\mu^{-1}(x)\bigr).
\end{equation}
Notice that $\pi$ is a discretization of $\pi_\infty$ in the sense of (\ref{pitilde}), because
\begin{equation}
{\pi}(x) := \pi_\infty(x_p) \,\mbox{ for }\, x_{p-1}< x\leq x_p.
\end{equation}
Therefore, by (\ref{Deltapitilde}),
\begin{align*}
\Delta(\pi, \pi_\infty) 
&\leq \kappa  \int_{0}^{1}\!  \bigl| \tr(\pi(x))-\tr(\pi_\infty(x)) \bigr| \, dx
\\
&=\kappa  \int_{0}^{1}\!  \bigl| \mu^{-1}(x)-\mu_\infty^{-1}(x)) \bigr| \, dx
\\
&= \kappa  \int_{0}^{1}\!  \bigl| \mu(x)-\mu_\infty (x)) \bigr| \, dx.
\end{align*}
In particular, we can choose the sequence in (\ref{xstrictagain}) to make $\Delta(\pi, \pi_\infty)$ as small as we want, while also making the distributions $\mu$ and $\mu_\infty$ as close as we want in $L^1$-norm.

As in Section \ref{Sec2label}, let $(v_\alpha)_{\alpha\in \Natural^r}$ be the weights of the Ruelle probability cascades corresponding to the parameters (\ref{xstrictagain}). Let $(\alpha^\ell)_{\ell\geq 1}$ be an i.i.d. sample from $\Natural^r$ according to these weights and, using the sequence of $q$'s in (\ref{xstrictagain}), define
\begin{equation}
T_{\ell,\ell'} = q_{\alpha^\ell\wedge \alpha^{\ell'}}.
\end{equation}
We already used in Section \ref{Sec6label} the fact that, by Theorem \ref{ThGGms}, the array $(\tr(R^{\infty}_{\ell,\ell'}))_{\ell,\ell'\geq 1}$ satisfies the classical Ghirlanda-Guerra identities. Therefore, Theorems 2.13 and 2.17 in \cite{SKmodel} imply that its distribution will be close to the distribution of the array $(T_{\ell,\ell'})_{\ell,\ell'\geq 1}$ when $\mu$ approximates $\mu_\infty.$ Consider the sequence  $\gamma_p:= \pi(x_p) = \Phi(q_p)$ for $0\leq p \leq r$, so that
\begin{equation}
0 =\gamma_0\leq  \ldots \leq \gamma_{r-1} \leq \gamma_r = \diag(d_1,\ldots,d_\kappa)
\label{zetafops}
\end{equation}
is a non-decreasing sequence in $\Gamma_\kappa$ (or $\Gamma_\kappa(d)$, satisfying the constraints in (\ref{LiftG})). Let
\begin{equation}
Q_{\ell,\ell'} = \Phi(T_{\ell,\ell'}) = \Phi \bigl(q_{\alpha^\ell\wedge \alpha^{\ell'}} \bigr).
\label{QLS}
\end{equation}
The fact that $\Phi$ is Lipschitz implies that the array $(Q_{\ell,\ell'})_{\ell,\ell'\geq 1}$ will be close in distribution to the array $(R^\infty_{\ell,\ell'})_{\ell,\ell'\geq 1}$.

Let us now consider Gaussian processes $Z^\alpha$ and $Y^\alpha$ indexed by $\alpha\in\Natural^r$ defined in Section \ref{Sec2label}, with the sequence (\ref{gammas}) now given by (\ref{zetafops}).  Consider a quantity similar to the quantity in (\ref{AS2repr}),
\begin{align}
f_M^{1,2}(\pi) &:= f_M^1(\pi) - f_M^2(\pi) 
\nonumber
\\
&:=
\frac{1}{M} \Bigl(
\e \log \sum_{\alpha\in\Natural^r} v_\alpha \sum_{\eps\in \Sigma_M(\delta^M)} \exp \beta \sum_{i\leq N} Z_{i}^\alpha(\eps_i)
-\e \log \sum_{\alpha\in\Natural^r} v_\alpha \exp \beta \sqrt{M} Y^\alpha
\Bigr).
\label{comp2}
\end{align}
If we compare the covariances in (\ref{Covz}) and (\ref{Covy}) with (\ref{CD}), Lemma \ref{LemCont2} implies that (\ref{comp2}) is the same continuous functional of the distribution of the array $(Q_{\ell,\ell'} )_{\ell,\ell'\geq 1}$ in (\ref{QLS}) as the quantity in (\ref{AS2repr}) is of the array $(R_{\ell,\ell'} )_{\ell,\ell'\geq 1}$ in (\ref{Rll}). Since both arrays, by construction, approximate in distribution the array $(R^{\infty}_{\ell,\ell'})_{\ell,\ell'\geq 1}$, we proved that the limit in (\ref{AS2repr}) equals to the limit of (\ref{comp2}) as $\mu \to \mu_\infty$ in $L_1$-norm. However, in this case $\Delta(\pi,\pi_\infty)\to 0$ and, by Lemma \ref{LemLC1} and (\ref{Lipschitz2}), this limit is just the extension $f_M^{1,2}(\pi_\infty)$ of the functional in (\ref{comp2}) to $\pi_\infty.$

To summarize, we showed that $f_M^{1,2}(\pi_\infty)$ gives the lower bound for the constrained free energy for any $M\geq 1$. Now we will let $M\to \infty$, but it is important to point out first that $\pi_\infty$ in the above construction depended on $M$, $\pi_\infty = \pi_\infty^M.$ Because the paths $\pi_\infty^M\in \Pi_d$ are monotone in $\Gamma_\kappa$, one can choose a convergent subsequence with respect to the metric $\Delta.$ Indeed, the diagonal elements of any $\pi\in\Pi_d$ are monotone functions, so one can choose an $L^1$-convergent subsequence of the diagonal elements first and then use that the sum $\pi_{k,k}+\pi_{k',k'}+2\pi_{k,k'}$ is also a monotone function, which allows to choose a convergent subsubsequence for the off-diagonal elements. Let $\pi^*$ be the limit of $\pi_\infty^M$ and let $\pi_\eps^*$ be a discretization of $\pi^*$ such that $\Delta(\pi^*, \pi_\eps^*)\leq \eps.$ By Lemma \ref{LemLC1} and (\ref{Lipschitz2}),
$$
\bigl| f_M^{1,2}(\pi_\infty^M) - f_M^{1,2}(\pi_\eps^*)\bigr| 
\leq
L \Delta(\pi_\infty^M, \pi_\eps^*)
\leq
L \bigl(\Delta(\pi_\infty^M, \pi^*) + \eps\bigr),
$$
and, therefore, $\limsup_{M\to\infty}|f_M^{1,2}(\pi_\infty^M) - f_M^{1,2}(\pi_\eps^*)| \leq L \eps.$
For discrete path $\pi^*_\eps,$ we can use Lemma \ref{Lem2}, which shows that the limit of the first term in (\ref{comp2}),
\begin{equation}
\lim_{M\to\infty} f_M^1(\pi_\eps^*) = \inf_{\lambda}\Bigl(
-\sum_{k\leq \kappa-1} \lambda_k d_k + \Phi(\lambda, d, \pi_\eps^*)
\Bigr).
\end{equation}
By (\ref{rearr2}), the second term is, actually, independent of $M$ and equals
\begin{equation}
f_M^2(\pi_\eps^*)
=
- \frac{\beta^2}{2}\sum_{k\leq \kappa} d_{k}^2 + \frac{\beta^2}{2}\int_{0}^1\! \|\pi_\eps^*(x)\|_{HS}^2\, dx.
\end{equation}
If we recall the functional $\PP(\lambda,d,\pi)$ defined in (\ref{PPldpi}), we proved that
$$
\liminf_{M\to\infty}  f_M^{1,2}(\pi_\infty^M)\geq \inf_{\lambda }\PP(\lambda,d,\pi_\eps^*) - L\eps
\geq  \inf_{\lambda,\pi\in \Pi_d}\PP(\lambda,d,\pi ) - L\eps.
$$
Letting $\eps\downarrow 0$ and maximizing over $d\in\DD$ finishes the proof of the lower bound.
\qed

\end{document}